\documentclass[11pt]{article}
\usepackage[latin1]{inputenc}
\usepackage{epsfig}
\usepackage{color}
\usepackage[british,english]{babel}
\usepackage{amsthm}
\usepackage{amsmath}
\usepackage{amsfonts}
\usepackage{amssymb}
\usepackage{graphicx}
\newtheorem{corollary}{Corollary}[section]
\newtheorem{definition}[corollary]{Definition}

\newtheorem{lemma}[corollary]{Lemma}

\newtheorem{remark}[corollary]{Remark}
\newtheorem{theorem}[corollary]{Theorem}
\newfont{\sBlackboard}{msbm10 scaled 900}

\newcommand{\mylabel}[1]{\label{#1}
            \ifx\undefined\stillediting
            \else \fbox{$#1$}\fi }
\newcommand{\BE}{\begin{equation}}

\newcommand{\EEQ}{\end{equation}}
\newcommand{\rfb}[1]{\mbox{\rm
   (\ref{#1})}\ifx\undefined\stillediting\else:\fbox{$#1$}\fi}

\newfont{\Blackboard}{msbm10 scaled 1200}

\newfont{\roma}{cmr10 scaled 1200}

\def\CC{\rm \hbox{C\kern-.56em\raise.4ex
         \hbox{$\scriptscriptstyle |$}\kern+0.5 em }}

\def\n{|\kern -.05cm{|}\kern -.05cm{|}}

\def \noame{\noalign{\medskip}}
\newcommand{\mm}    {{\hbox{\hskip 0.5pt}}}

\newcommand{\bluff} {{\hbox{\raise 15pt \hbox{\mm}}}}

\newcommand{\ep}   {\varepsilon}
\usepackage{fancyhdr}

\makeatletter
\def\section{\@startsection {section}{1}{\z@}{-3.5ex plus -1ex minus
    -.2ex}{2.3ex plus .2ex}{\large\bf}}
\makeatother
\def\be{\begin{equation}}
\def\ee{\end{equation}}

\date{ }
\begin{document}
\thispagestyle{empty}
\title{\Large \bf    Asymptotic analysis of the Navier-Stokes equations in a thin domain\\ with power law slip boundary conditions}\maketitle
\vspace{-2cm}
\begin{center}
Mar\'ia ANGUIANO\footnote{Departamento de An\'alisis Matem\'atico. Facultad de Matem\'aticas. Universidad de Sevilla. 41012-Sevilla (Spain) anguiano@us.es} and Francisco Javier SU\'AREZ-GRAU\footnote{Departamento de Ecuaciones Diferenciales y An\'alisis Num\'erico. Facultad de Matem\'aticas. Universidad de Sevilla. 41012-Sevilla (Spain) fjsgrau@us.es}
 \end{center}

 \renewcommand{\abstractname} {\bf Abstract}
\begin{abstract} 
This theoretical study deals with the Navier-Stokes equations posed in a 3D thin domain with thickness $0<\ep\ll 1$, assuming power law slip boundary conditions, with an anisotropic tensor, on the bottom. This condition,  introduced in (Djoko {\it et al.} {\it Comput. Math. Appl.}  128 (2022) 198--213), represents a generalization  of the Navier slip boundary condition. The goal  is to study the influence of the power law slip boundary conditions with an anisotropic tensor of order $\ep^{\gamma\over s}$, with $\gamma\in \mathbb{R}$ and  flow index $1<s<2$, on the behavior of the fluid with thickness $\ep$ by using asymptotic analysis when $\ep\to 0$, depending on the values of $\gamma$. As a result, we deduce the existence of a critical value of $\gamma$ given by $\gamma_s^*=3-2s$ and so, three different limit boundary conditions are derived. The critical case  $\gamma=\gamma_s^*$ corresponds to a  limit condition of type power law slip. The supercritical case $\gamma>\gamma_s^*$ corresponds to a limit boundary condition of type perfect slip. The subcritical case  $\gamma<\gamma_s^*$ corresponds to a limit  boundary condition of type no-slip.

 \end{abstract}
\bigskip\noindent

\noindent {\small \bf AMS classification numbers:}  35Q35, 76A20, 76A05, 76M50.
\\
 

\section{Introduction}
The stationary Navier-Stokes equations in a domain $\Omega$ reads as follows
\begin{equation}\label{Navierequation}-2\nu\,{\rm div}(\mathbb{D}[u])+(u\cdot \nabla)u+\nabla p=f\quad\hbox{and}\quad {\rm div}(u)=0,
\end{equation}
where $u$ denotes the velocity field,  $\mathbb{D}[u]={1\over 2}(Du+(Du)^T)$   the deformation tensor associated with the velocity field $u$, $p$ the scalar pressure, $f$ the external forces and $\nu>0$ the viscosity.
Concerning the boundary conditions, it is commonly accepted that viscous fluids adhere to surfaces, and so the no-slip condition at the surfaces of a domain, given  by 
$$u=0\quad\hbox{on}\quad\partial\Omega,$$ 
is widely used.   Under suitable regularity conditions on the domain and $f$, this problem is well studied ma\-thematically, see for instance Boyer \& Fabrie \cite{Boyer},  Galdi \cite{Galdi} or Temam \cite{Temam}.  However, this condition does not seem always valid physically, indeed some fluids melt and solutions slip against the surface. Also, sometimes the no-slip condition is not good enough because it is not possible to describe the behavior of the fluid near the boundary. Therefore, it is necessary to introduce other type of boundary conditions to describe this behavior. In this sense, Navier \cite{Navier} proposed the Navier slip boundary conditions in which it is assumed a thin layer of a fluid near of the boundary and the tangential component of the strain tensor should be proportional to the tangential component of the fluid velocity on a part of the boundary $\Gamma\subset\partial\Omega$, that is
\begin{equation}\label{Navierslip}
2\nu[\mathbb{D}[u]{\rm n}]_\tau=-\lambda[u]_\tau,\quad u\cdot {\rm n}=0,\quad\hbox{on }\Gamma,
\end{equation}
where ${\rm n}$ denotes the outside unitary normal vector to $\Omega$ on $\Gamma$, $\lambda> 0$ is the friction coefficient and  the subscript $\tau$ denotes the orthogonal projection on the tangent space of $\Gamma$, i.e.  $[u]_\tau=u-(u\cdot {\rm n}){\rm n}$.  Problem (\ref{Navierequation}) with Navier slip boundary conditions (\ref{Navierslip})  has been studied by many authors in different contexts, see for example Amrouche  \& Rejaiba \cite{Amrouche}, Clopeau {\it et al.} \cite{Clupeau} and  Solonnikov  \& ${\rm \check{S}}$${\rm \check{c}}$adilov \cite{Solonnikov}.  Notice that depending on the value of $\lambda$ in (\ref{Navierslip}), we shall consider the following type of boundary conditions:
\begin{itemize}
\item  Perfect slip when $\lambda=0$, i.e. 
\begin{equation}\label{pure_slip}
2\nu[\mathbb{D}[u]{\rm n}]_\tau=0,\quad u\cdot {\rm n}=0,\quad\hbox{on }\Gamma,
\end{equation}
\item Partial slip when $\lambda\in (0,+\infty)$,  
\begin{equation}
2\nu[\mathbb{D}[u]{\rm n}]_\tau=-\lambda [u]_\tau,\quad u\cdot {\rm n}=0,\quad\hbox{on }\Gamma,
\end{equation}
\item No-slip when $\lambda=+\infty$, i.e. 
$$[u]_\tau=0, \quad u\cdot {\rm n}=0,\quad\hbox{on }\Gamma,$$
which implies $u=0$ on $\Gamma$.
\end{itemize}

Related to this, we refer to Acevedo {\it et al.} \cite{Acevedo} for the study the limiting behavior of the solution $(u_\lambda, p_\lambda)$ of problem  (\ref{Navierequation}) with Navier slip boundary conditions (\ref{Navierslip}), when the friction coefficient $\lambda$ goes to $0$ or $\infty$. In fact, they proved that $(u_\lambda, p_\lambda)$ weakly converges to $(u_0, p_0)$ when $\lambda\to 0$ in suitable Sobolev spaces, where $(u_0, p_0)$ is the solution of the Navier-Stokes system with Navier slip boundary conditions corresponding with $\lambda=0$.  Also, it holds that $(u_\lambda, p_\lambda)$ weakly converges to $(u_\infty, p_\infty)$ when $\lambda\to \infty$, where $(u_\infty, p_\infty)$ is the solution of the Navier-Stokes system with no-slip boundary conditions, i.e. the Navier boundary conditions corresponding with $\lambda=+\infty$.

In this work,  we are interested in a generalization of the Navier slip condition recently introduced by Djoko  {\it et al.} \cite{Djoko} (see also Aldbaissy {\it et al.} \cite{Aldbaissy, Aldbaissy2} and Djoko {\it et al.} \cite{Djoko2}), which  arises when the contact surface is lubricated with a thin layer of a non-Newtonian fluid. This condition is called power law slip boundary condition and reads as follows
\begin{equation}\label{Powerslip}
2\nu[\mathbb{D}[u]{\rm n}]_\tau=-|K [u]_\tau|^{s-2}K^2[u]_\tau,\quad u\cdot {\rm n}=0,\quad\hbox{on }\Gamma,
\end{equation}
where $|v|^2=v\cdot v$ is the Euclidean norm. We observe that in this condition, the tangential shear is a power law function of the tangential velocity, where $K\in\mathbb{R}^{2\times 2}$ is an anisotropic tensor, assumed to be uniformly positive definite, symmetric and bounded, and  $s$ is a real, strictly positive number representing the flow behavior index.  

 The boundary condition  (\ref{Powerslip}) represents a generalization of the Navier slip boundary condition (\ref{Navierslip}), since for $s=2$ and $K=\lambda^{1\over 2}I$ with $\lambda>0$, then the power slip boundary condition (\ref{Powerslip}) reduces to Navier slip boundary condition (\ref{Navierslip}).    We also mentioned that the power slip boundary condition (\ref{Powerslip}) is present in the context of laminar flows of Newtonian liquids (e.g. water) over complex surfaces, also when a rough or structured boundary surface is anisotropic, e.g. when it has rows of riblets, pillars or periodic patterns, the effective slip condition is anisotropic, i.e., direction dependent. When the surface is heterogeneous, the effective slip is also position-dependent. This can occur, for example, when the boundary has a varying degree of roughness or when the boundary is a smooth surface with a varying hydrophobic/hydrophilic composition.  For instance, we refer to the derivation of effective slip boundary conditions coming from rough boundaries in Bonnivard {\it et al.} \cite{Bonnivard}, Bonnivard \& Su\'arez-Grau \cite{Bonnivard_SG1, Bonnivard_SG2}, Bucur \cite{BucurNN}, Bucur {\it et al.} \cite{Bucur},  Dalibard \& G\'erard-Varet \cite{Dalibard} and  Su\'arez-Grau \cite{SG1, SG2}.  The existence of solutions of the Stokes and Navier-Stokes equations with power law boundary conditions (\ref{Powerslip})  on a part of the boundary was studied in \cite{Djoko}  for $1<s<2$, which corresponds to the tangential shear thinning. In the  case $s>2$ the existence of solutions is not proven, it was not able to prove a inf-sup conditon, which is the key point to obtain the pressure. In the case $s=2$, repeating the classical proof of the existence of solution of  the Stokes and Navier-Stokes problem with homogeneous Dirichlet conditions (see for instance   \cite[Theorem 2.3]{Acevedo}, \cite[Theorem 7.1]{Lions} and \cite[Theorem 10.1]{Temam})  gives the existence of solution of the Stokes and Navier-Stokes equations with Navier slip boundary conditions.

Our main interest in this paper is to study a lubrication problem corresponding to the asymptotic influence of the power law boundary condition (\ref{Powerslip}), imposed on a part of the boundary  with $1<s<2$,  on the  behavior of the Navier-Stokes equations  through a thin domain $\Omega^\ep$, where the small parameter $0<\ep\ll 1$ represents the thickness of the domain. More precisely, we consider the following $3D$ thin domain  (see Figure \ref{fig:thin domain})
$$\Omega^\varepsilon=\{(x_1,x_2, x_3)\in\mathbb{R}^2\times \mathbb{R}\,:\, (x_1,x_2)\in \omega,\quad 0<x_3<\ep h(x_1,x_2)\},$$
where $\omega$ is a smooth, connected open set of $\mathbb{R}^2$ and $h$ is a smooth and positive function (see Section \ref{S1} for more details). To study the influence of the slip boundary conditions on the behavior of the Navier-Stokes equations in the thin domain $\Omega^\ep$ (the subscript $\ep$ is added to the unknowns to stress the dependence of the solution on the small parameter)
$$
\left\{\begin{array}{rl}
- \nu\,\Delta u_\ep+(u_\ep\cdot \nabla )u_\ep  +\nabla p_\ep= f_\ep &\hbox{in } \Omega^\ep,\\
\noame
{\rm div}(   u_\varepsilon)=0& \hbox{in } \Omega^\ep,
\end{array}\right.
$$ 
 where we consider the case of the power law slip boundary conditions (\ref{Powerslip}) on $\Gamma_0$ with an anisotropic tensor $K_\ep$, depending on $\varepsilon$, of the form
$$K_\ep=\ep^{\gamma\over s}K,\quad\hbox{with } 1<s<2\quad\hbox{and}\quad \gamma\in\mathbb{R},$$
where $K$ is assumed to be uniformly positive definite, symmetric and bounded. Thus, the power slip boundary condition, 
 where the tangential shear is a power law function of the tangential velocity with a coefficient depending on $\varepsilon$, is given by
\begin{equation}\label{power_slip_ep}2[\mathbb{D}[u_\ep]{\rm n}]_\tau=-\ep^{\gamma} |K [u_\ep]_\tau|^{s-2}K^2[u_\ep]_\tau,\quad u_\ep\cdot {\rm n}=0,\quad\hbox{on  }\Gamma_0,
\end{equation}
and no-slip condition on the rest of the boundary, i.e.
$$u_\ep=0\quad   \hbox{on }  \partial\Omega^\ep\setminus \Gamma_0.$$
 After the homogenization process (under assumptions given in Section \ref{S1}) when $\ep\to0$ depending on the value of $\gamma$, we deduce (see Theorem \ref{thm1_critical_rig}) that the limit velocity $\widetilde u=(\widetilde u',0)$ and limit pressure $\widetilde p$ satisfies the reduced 2D-Stokes system   
 \begin{equation}\label{limit_crit_rig_intro}
\left\{\begin{array}{rl}
\displaystyle - \nu\partial_{z_3}^2\widetilde u' (z)=  f'(z')-\nabla_{z'}\widetilde p(z') & \hbox{in } \Omega=\{z\in\mathbb{R}^3\,: \, z'\in \omega,\  0<z_3<h(z')\},\\
\noame
\displaystyle{\rm div}_{z'}\left(\int_0^{h(z')}\widetilde u'(z)\,dz_3\right)=0& \hbox{in } \omega,\\
\noame
\displaystyle \left(\int_0^{h(z')}\widetilde u'(z)\,dz_3\right)\cdot {\rm n} =0& \hbox{on } \partial\omega,\\
\noame
\displaystyle \widetilde u'=0&\hbox{on }\Gamma_1=\omega\times \{h(z')\},
\end{array} \right.
\end{equation}
 where $\widetilde u'=(u_1, u_2)$ and $f'=(f_1, f_2)$. Moreover, we prove the existence of a critical value for $\gamma$ given by
 \begin{equation}\label{gammastar}
\gamma_s^*=3-2s,\quad \hbox{with}\quad 1<s<2,
\end{equation}
which let us derive three different boundary conditions for $\widetilde u'$ on the bottom $\Gamma_0$:  
\begin{itemize}
\item If $\gamma=\gamma_s^*$, then the effective boundary condition on $\Gamma_0$ is a power slip boundary condition with anisotropic  tensor $K$, i.e.
 $$-\nu\partial_{z_3}\widetilde u'=-|K\widetilde u'|^{s-2}K^2\widetilde u'\quad\hbox{on }\Gamma_0.$$ 
 Thus, to take into account the anisotropy, then $K_\ep$ has to be of order $\mathcal{O}(\ep^{\gamma_{s^*}\over s})$.
\item If $\gamma>\gamma_s^*$, then the effective boundary condition on $\Gamma_0$ is the perfect slip boundary condition, i.e. 
$$-\nu\partial_{z_3}\widetilde u'=0 \quad\hbox{on }\Gamma_0.$$ 
 This means that for an anisotropy tensor of order smaller than $\mathcal{O}(\ep^{\gamma_s^*\over s})$, then the fluid does not take into account anisotropy and slides perfectly.
\item If $\gamma<\gamma_s^*$, then the effective boundary condition on $\Gamma_0$ is the no-slip condition, i.e.
$$\widetilde u'=0\quad\hbox{on }\Gamma_0.$$  This means that for an anisotropy tensor of order greater than $\mathcal{O}(\ep^{\gamma_s^*\over s})$, then the anisotropy is so strong that the fluid is stopped on the boundary.
\end{itemize}

Observe that for $s=2$ and $K_\ep=\ep^{\gamma\over 2}\lambda^{1\over 2} I$, with $\lambda>0$, where the power slip condition (\ref{power_slip_ep}) reduces to the Navier slip condition with friction parameter $\lambda\ep^\gamma$,  it holds that the critical value $\gamma_2^*=-1$, which is the critical value for the case of Navier slip boundary condition. Namely,    to take into account the friction coefficient $\lambda$ in the effective boundary condition, i.e.
$$-\nu\partial_{z_3}\widetilde u'=-\lambda \widetilde u'\quad\hbox{on }\Gamma_0,$$
the original friction coefficient has to be of order $\mathcal{O}(\ep^{-1})$. If the original friction coefficient $\lambda \ep^\gamma$ is of order smaller than $\mathcal{O}(\ep^{-1})$, then the fluid behaves on the boundary as if there were no friction (perfect slippage), i.e.
$$-\nu\partial_{z_3}\widetilde u'=0\quad\hbox{on }\Gamma_0.$$
Finally, if the friction coefficient $\lambda \ep^\gamma$ is of order greater than $\mathcal{O}(\ep^{-1})$, then the friction coefficient is so strong that  the fluid is stopped on the boundary (no-slip condition), i.e. 
$$\widetilde u'=0\quad\hbox{on }\Gamma_0.$$

  To prove these results, we first use the multiscale expansion method, which is a formal but powerful tool to analyse homogenization problems, see for instance the application of this method in Bayada \& Chambat  \cite{Bayada} and  Mikeli\'c \cite{Mikelic}. Next, once the results have been understood, we rigorously justify them by means of the derivation of  {\it a priori} estimates and  some compactness results.

As far as the authors know, the flow of a Newtonian fluids with power law slip boundary conditions has not been
yet considered in the above described lubrication framework, which represents the main novelty of the paper.    We observe that the obtained findings are amenable for the numerical simulations with a considerable simplification with respect to the original problem (which is computationally more expensive), since the effective system (\ref{limit_crit_rig_intro}) is a two dimensional ordinary differential system with respect to $z_3$. Therefore,  we believe that it could prove useful in the engineering practice as well.

The paper is structured as follows. In Section \ref{S1}, we introduce the statement of the problem. In Section \ref{sec:formal}, we consider the formal derivation, and in Section \ref{sec:Rigorous} we will rigorously justify the results. We finish the paper with a section of references.

\section{Formulation of the problem and preliminaries}\label{S1}
In this section, we first define the thin domain and some sets necessary to study the asymptotic behavior of the solutions. Next, we introduce the problem considered in the thin domain and also, the rescaled problem posed in the domain of fixed height, together with the respective weak variational formulations.

\paragraph{The domain and some notation.} 
Along this paper, the points $x\in\mathbb{R}^3$ will be decomposed as $x=(x',x_3)$ with $x'\in\mathbb{R}^2$, $x_3\in\mathbb{R}$. We also use the notation $x'$ to denote a generic vector of $\mathbb{R}^2$.

\begin{figure}[h!]
\begin{center}
\includegraphics[width=8cm]{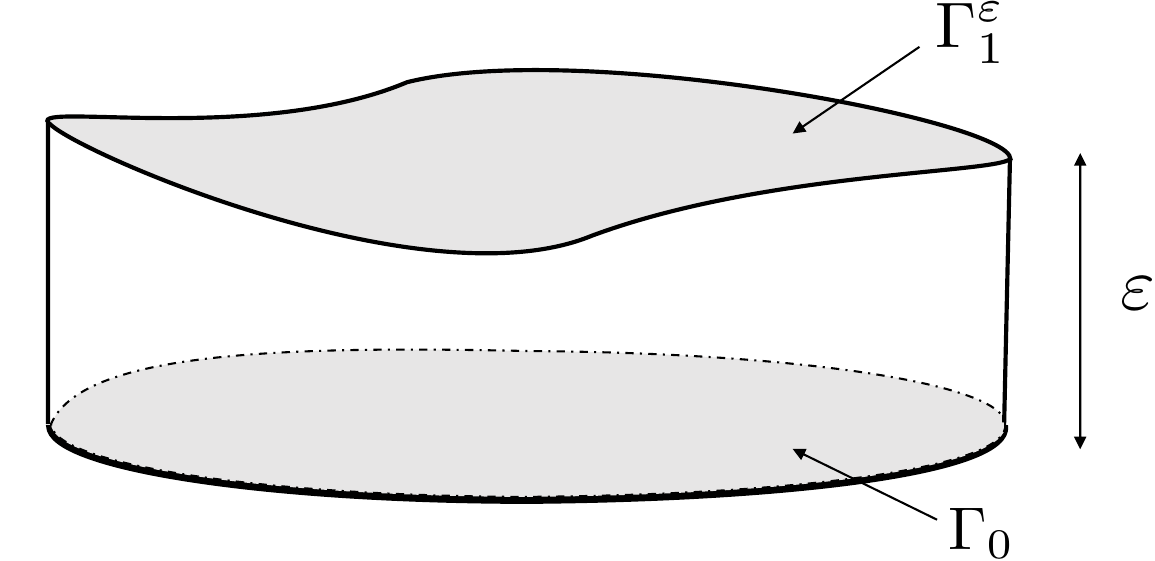}
\end{center}
\caption{Thin domain $\Omega^\ep$, top boundary $\Gamma^\ep_1$ and bottom boundary $\Gamma_0$}
\label{fig:thin domain}
\end{figure}

\noindent We consider  $\omega$ as an open, smooth, bounded and
connected set of $\mathbb{R}^2$, and  a 3D thin domain given by
$$\Omega^\varepsilon=\{(x', x_3)\in\mathbb{R}^2\times \mathbb{R}\,:\, x'\in \omega,\quad 0<x_3<h_\ep(x')\},$$
Here, the function $h_\ep(x')=\ep h(x')$ represents the real gap between the two surfaces. The small para\-meter $\ep$ is related to the film thickness.  Function $h$ is positive and smooth  $C^1$ bounded function defined for $x'$.  We define the bottom, top and lateral boundaries of $\Omega^\ep $ as follows (see Figure \ref{fig:thin domain})
$$\Gamma_0=\omega\times \{0\},\quad \Gamma^\ep _1=\left\{(x',x_3)\in\mathbb{R}^3\,:\, x'\in \omega,\ x_3= h_\ep(x')\right\},\quad \Gamma_\ell^\ep=\partial\Omega^\ep\setminus (\Gamma_0\cup \Gamma_1^\ep).$$

Let us now introduce some notation which will be useful in the following. For a vectorial function $ \varphi=(  \varphi',  \varphi_3)$ and a scalar function $ \phi$, we introduce the operators $\Delta$,  ${\rm div}$, $D$ and $\nabla$    by
\begin{equation}\label{def1}\begin{array}{c}
\Delta\varphi=\Delta_{x'}\varphi+\partial_{x_3}^2 \varphi, \quad {\rm div} (\widetilde\varphi)={\rm div}_{x'}( \varphi')+ \partial_{x_3} \varphi_3,\\
\noame
\displaystyle
(D \varphi)_{ij}=\partial_{x_j} \varphi_i\ \hbox{for }i=1,2,3,\ j=1,2, 3, 
\\
\noame
\displaystyle
\nabla  \phi =(\nabla_{x'}  \phi, \partial_{x_3} \phi)^t.
\end{array}
\end{equation}
\noindent We denote by $O_\ep$ a generic real sequence which tends to zero with $\ep$ and can change from line to
line. We denote by $C$ a generic positive constant which can change from line to line.

\paragraph{The problem and the rescaling.} 
 As stated in the introduction, we consider the 3D stationary Navier-Stokes equations by setting 
$$u_\ep=(u'_\ep(x), u_{3,\ep}(x)),\quad p_\ep=p_\ep(x),$$
at a point $x\in \Omega^\ep$, which is given by
\begin{equation}\label{system_1introad}
\left\{\begin{array}{rl}
- \nu\,\Delta u_\ep+(u_\ep\cdot \nabla )u_\ep  +\nabla p_\ep= f_\ep &\hbox{in } \Omega^\ep,\\
\noame
{\rm div}(   u_\varepsilon)=0& \hbox{in } \Omega^\ep,\\
\noame
  u_\ep=0&   \hbox{on }  \Gamma_1^\ep\cup \Gamma_\ell^\ep,
\end{array}\right.
\end{equation}
where $\nu>0$, and the power law slip boundary conditions  prescribed on $\Gamma_0$ given by
\begin{equation}\label{bc0_system}
-\nu\partial_{x_3} u_\ep'=-  \ep^\gamma |K u_\ep'|^{s-2}K^2 u_\ep' ,\quad  u_{\ep,3}=0,\quad  \hbox{on }   \Gamma_0,
\end{equation}
 where $1<s<2$ and $\gamma\in \mathbb{R}$.

\begin{remark} Here, we have assumed the following:
\begin{itemize} 
\item Since ${\rm div}(u_\ep)=0$, it holds $2\,{\rm div}(\mathbb{D}[u_\ep])=\Delta u_\ep$. Then, it also holds $2\nu[\mathbb{D}[u_\ep]\,{\rm n}]_\tau=\nu[Du_\ep\,{\rm n}]_\tau$.

\item Since the bottom boundary $\Gamma_0$ is flat, the outside normal vector ${\rm n}=-{\rm e}_3$, where $\{{\rm e}_i\}_{i=1}^3$ is the canonical basis in $\mathbb{R}^3$, and so $u_\ep\cdot {\rm n}=0$ implies $u_{\ep,3}=0$.  Also,  the orthogonal projection of a function $u_\ep$ on the tangent space of $\Gamma_0$ is $[u_\ep]_\tau=u'_\ep$, where $u'_\ep=(u_{\ep,1}, u_{\ep,2})$, and then, $\nu[Du_\ep\,{\rm n}]_\tau=-\nu\partial_{x_3} u'_\ep$.

\item Due to the thickness of the domain, it is usual to assume that the vertical
components of the external forces can be neglected and, moreover, the forces can be
considered independent of the vertical variable. Thus, for sake of simplicity, given $f'=(f_1, f_2)\in L^2(\omega)^2$, along the paper we consider the following type of external forces $f_\ep $ (see for instance \cite{SG1, SG2}): 
$$f_\ep=(f'(x'),0)^t.$$
\end{itemize}

\end{remark}

\begin{definition}\label{def:weakform} For $\ep>0$,  we say that  $(u_\ep,  p_\ep)$ defined on $\Omega^\ep$ is a weak solution of problem (\ref{system_1introad})--(\ref{bc0_system}) if and only if the functions  $(u_\ep,p_\ep)\in \mathbb{V}(\Omega^\ep) \times L^2_0(\Omega^\ep)$,  where the corresponding functional space for velocity is 
$$ \mathbb{V}(\Omega^\ep)=\left\{\varphi\in H^1(\Omega^\ep)^3\,:\,  \varphi=0 \hbox{ on }\partial\Omega^\ep\setminus \Gamma_0,\  \ \varphi_3=0\quad\hbox{on }\Gamma_0\right\},$$
and the space for pressure $L^{2}_0$ is the space of functions of $L^2$ with zero mean value, and satisfy 
\begin{equation}\label{form_var_vel}
\begin{array}{l}
\displaystyle
\nu \int_{\Omega^\ep}Du_\ep:D\varphi\,dx+\int_{\Omega^\ep}(u_\ep\cdot \nabla )u_\ep\cdot \varphi\,dx+\ep^\gamma \int_{\Gamma_0}|Ku_\ep'|^{s-2}K  u_\ep'\cdot K \varphi'\,d\sigma-\int_{\Omega^\ep}p_\ep\,{\rm div}(\varphi)\,dx\\
\noame
\displaystyle =\int_{\Omega^\ep}f'\cdot \varphi'\,dx,\quad
  \forall\,\varphi\in\mathbb{V}(\Omega^\ep),
 \end{array}
\end{equation}
and
\begin{equation}\label{formdiv}
\int_{\Omega^\ep}{\rm div}(u_\ep)\, \psi\,dx=0,\quad \forall\,\psi\in L^2(\Omega^\ep).
\end{equation}
\end{definition}
  
\begin{remark} Under previous assumptions, for every $\ep>0$,  we have that  reference \cite[Proposition 2.2]{Djoko} gives the existence of at least one weak solution $(u_\ep, p_\ep)$ of problem (\ref{system_1introad})--(\ref{bc0_system}).
\end{remark}

\noindent The objetive of this paper is to study the asymptotic problems   for the behavior of  the sequence of solutions $(u_\ep,  p_\ep)$ of previous problems, when $\ep$ tends to zero depending on the value of $\gamma$.  To do that, we introduce a classical change of variables in thin domains,  the dilatation
\begin{equation}\label{dilatation}
z'=x',\quad z_3=\ep^{-1} x_3.
\end{equation}
This change transforms $\Omega^\ep$ into a fixed domain $\Omega$, defined by
\begin{equation}\label{Omegatildeep}
 \Omega=\left\{(z',z_3)\in \mathbb{R}^2\times \mathbb{R}\,:\, z'\in \omega,\ 0<z_3<h(z')\right\}.
\end{equation}
The boundary of $\Omega$ is denoted by $\partial\Omega$, where the top and lateral boundaries of the rescaled domain $\Omega$ is defined by
\begin{equation}\label{OmegatildeepGamma}\Gamma_1=\left\{(z',z_3)\in \mathbb{R}^2\times \mathbb{R}\,:\, z'\in \omega,\ z_3=h(z')\right\}, \quad \Gamma_\ell=\partial\Omega\setminus (\Gamma_0\cup \Gamma_1).
\end{equation}
Accordingly, we define the functions $\widetilde u_\ep$ and $\widetilde p_\ep$ by
\begin{equation}\label{tildefunctions}
\widetilde u_\ep (z)=u_\ep (z',\ep z_3),\quad \widetilde p_\ep (z)=p_\ep(z',\ep z_3)\quad \hbox{a.e. }z\in \Omega.
\end{equation}
Observe that according to the assumption on $f_\ep$, it holds $\widetilde f_\ep(z)=(f'(z'),0)$ a.e. $z\in \Omega$.\\

Let us now introduce some notation which will be useful in the following. For a vectorial function $\widetilde \varphi=(\widetilde \varphi',\widetilde \varphi_3)$ and a scalar function $\widetilde \phi$ obtained respectively from functions  $\varphi$ and $\phi$ by using the change of variables (\ref{dilatation}), we introduce the operators $\Delta_\ep$,  ${\rm div}_\ep$, $D_\ep$ and $\nabla_\ep$ by
\begin{equation}\label{def_divrot}\begin{array}{c}
\Delta_\ep\widetilde\varphi=\Delta_{z'}\widetilde\varphi+\ep^{-2}\partial_{z_3}^2\widetilde\varphi, \quad {\rm div}_\ep (\widetilde\varphi)={\rm div}_{z'}\widetilde\varphi'+\ep^{-1}\partial_{z_3}\widetilde\varphi_3,\\
\noame
\displaystyle
(D_\ep\widetilde\varphi)_{ij}=\partial_{z_j}\widetilde\varphi_i\quad \hbox{for }i=1,2,3,\ j=1,2,\quad (D_\ep\widetilde\varphi)_{i3}=\ep^{-1}\partial_{z_3}\widetilde\varphi_i\quad \hbox{for }i=1,2,3,
\\
\noame
\displaystyle
\nabla_\ep\widetilde\phi =(\nabla_{z'}\widetilde \phi,\ep^{-1}\partial_{z_3}\widetilde\phi).
\end{array}
\end{equation}

\noindent Thus, using the change of variables (\ref{dilatation}), the system (\ref{system_1introad})--(\ref{bc0_system}) can be rewritten as
\begin{equation}\label{system_1introtilde}
\left\{\begin{array}{rl}
- \nu\Delta_\ep \widetilde u_\ep+ (\widetilde u_\ep\cdot \nabla_\ep)\widetilde u_\ep +\nabla_\ep\widetilde  p_\ep= \widetilde f_\ep &\hbox{in } \Omega,\\
\noame
{\rm div}_\ep(  \widetilde  u_\varepsilon)=0& \hbox{in } \Omega,\\
\noame
\displaystyle  \widetilde u_\ep=0&   \hbox{on }  \Gamma_1 \cup \Gamma_\ell,\\
\end{array}\right.
\end{equation}
with power law  slip  boundary conditions
\begin{equation}\label{bc1_system_tilde_0}
 \begin{array}{l}
 \displaystyle -{\nu\over \ep}\partial_{z_3} \widetilde u_\ep'=-  \ep^\gamma |K \widetilde u_\ep'|^{s-2}K^2 \widetilde u_\ep',\quad \widetilde u_{\ep, 3}=0\quad     \hbox{on }   \Gamma_0.
 \end{array} 
\end{equation}

\noindent According to the change of variables (\ref{dilatation}) applied to the weak variational formulations given in Definition \ref{def:weakform}, then, for $\ep>0$,   a rescaled weak solution $(\widetilde u_\ep, \widetilde p_\ep)\in \mathbb{V}(\Omega)\times L^2_0(\Omega)$, where 
the corresponding functional space for velocity is 
$$ \mathbb{V}(\Omega)=\left\{\varphi\in H^1(\Omega)^3\,:\,  \varphi=0 \hbox{ on }\partial\Omega\setminus \Gamma_0,\  \ \varphi_3=0\quad\hbox{on }\Gamma_0\right\},$$ satisfies 
\begin{equation}\label{form_var_vel_tilde}
\begin{array}{l}
\displaystyle
 \nu \int_{\Omega}D_\ep \widetilde u_\ep :D_\ep\widetilde \varphi\,dz+\int_{\widetilde \Omega^\ep}(\widetilde u_\ep\cdot \nabla_\ep ) \widetilde u_\ep\cdot \widetilde \varphi\,dz+\ep^{\gamma-1} \int_{\Gamma_0}|K\widetilde u_\ep'|^{s-2} K\widetilde u_\ep'\cdot K\widetilde \varphi'\,d\sigma
 \displaystyle -\int_{\Omega}\widetilde p_\ep\,{\rm div}_\ep(\widetilde \varphi)\,dz\\
 \noame \displaystyle=\int_{\Omega}f'\cdot \widetilde \varphi'\,dz,
 \end{array}
\end{equation}
and
\begin{equation}\label{formdiv_tilde}
\int_{\Omega}{\rm div}_\varepsilon(\widetilde u_\ep) \,\widetilde \psi\,dz=0,
\end{equation}
  for every  $\widetilde \varphi\in \mathbb{V}(\Omega)$ and $\widetilde\psi\in L^2(\Omega)$  obtained from $(\varphi, \psi)$ by the change of variables (\ref{dilatation}).

Now, the goal is to study the asymptotic behavior of a sequence of solution $(\widetilde u_\ep, \widetilde p_\ep)$ of problem (\ref{system_1introtilde})--(\ref{bc1_system_tilde_0}). In the next section, we will study the asymptotic analysis in a formal way, and in Section \ref{sec:Rigorous}, we develop the rigorous analysis.

\section{The formal asymptotic expansion}\label{sec:formal}
In this section, we apply the asymptotic expansion method (see for instance \cite{Bayada, Mikelic}) to the Navier-Stokes equations (\ref{system_1introtilde}) with power slip boundary conditions (\ref{bc1_system_tilde_0}). We will derive a reduced Stokes system with no-slip condition on the top boundary $\Gamma_1$ and different boundary conditions on $\Gamma_0$ depending on the value of $\gamma$. The idea is to assume an expansion in $\ep$ of the solution $(\widetilde u_\ep, \widetilde p_\ep)$ given by
\begin{equation}\label{expansionbeta}\widetilde u_\ep(z)=\ep^\beta\Big(v^0(z)+ \ep v^1(z)+\ep^2 v^2(z)+\cdots\Big) ,\quad \widetilde p_\ep(z)=p^0(z)+\ep p^1(z) +\ep^2 p^2(z)+\ep^3 p^3(z)+\cdots
\end{equation}
To determine the effective problem given by functions  $(v^0, p^0)$, the expansion (\ref{expansionbeta}) is plugged into the PDE, we identify the various powers of $\ep$ and we obtain a cascade of equations from which we retain only the leading ones that constitute the effective problem. 

We remark that the value $\beta$ in the expansion (\ref{expansionbeta}) will be determined in the next section, by deriving {\it a priori} estimates for $(\widetilde u_\ep, \widetilde p_\ep)$ (see Lemmas \ref{lem:estimvelmic} and \ref{lem:estimpress}).  Moreover, the effective problem will be justified by corresponding compactness results (see Lemma \ref{lemmaconv_crit} and Theorem \ref{thm1_critical_rig}).

\begin{theorem} \label{thm:formal_NS}
Assume $1<s<2$ and define $\gamma^*_s=3-2s$. Consider $(\widetilde u_\ep, \widetilde p_\ep)$ a sequence of solutions of problem  (\ref{system_1introtilde})--(\ref{bc1_system_tilde_0}). Assuming the asymptotic expansion of the unknown $(\widetilde u_\ep, \widetilde p_\ep)$ in the following form
\begin{equation}\label{formal_exp_Nav}
\widetilde u_\ep(z)=\ep^2 v^0(z)+ \ep^3 v^1(z)+\ep^4v^2(z)+\cdots ,\quad \widetilde p_\ep(z)=p^0(z)+\ep p^1(z) +\ep^2 p^2(z)+\ep^3 p^3(z)+\cdots
\end{equation}
for a.e. $z\in \Omega$, where 
$$v^i=( \bar v^i, v_3^i)\quad \hbox{with}\quad \bar v^i=(v_1^i,v_2^i),\quad i=0, 1, \ldots,$$ 
we deduce that the main order pair of functions $(v^0, p^0)$, with $v_3^0\equiv 0$ and $p^0=p^0(z')$, satisfies the following effective reduced Stokes  problem
\begin{equation}\label{limit_crit_Nav}
\left\{\begin{array}{rl}
\displaystyle - \nu\partial_{z_3}^2 \bar v^0(z)=  f'(z')-\nabla_{x'}p^0(z')&  \hbox{in } \Omega,\\
\noame
\displaystyle {\rm div}_{z'}\left(\int_0^{h(z')}\bar v^0(z)\,dz_3\right)=0& \hbox{in } \omega,\\
\noame
\bar v^0=0&\hbox{ on }\Gamma_1,
\end{array} \right.
\end{equation}
Moreover, $\bar v^0$ satisfies the following boundary condition on the bottom $\Gamma_0$ depending on the value of $\gamma$:
\begin{itemize}
\item If $\gamma=\gamma_s^*$, then it holds a power law slip boundary condition
\begin{equation}\label{crit_bc_limit_formal}
 -\nu\partial_{z_3}\bar v^0=-|K\bar v^0|^{s-2}K^2\bar v^0\quad \hbox{ on }\Gamma_0.
\end{equation}

\item $\gamma>\gamma_s^*$, then it holds a perfect slip boundary condition
\begin{equation}\label{super_bc_limit_formal}
 -\nu\partial_{z_3}\bar v^0=0\quad \hbox{ on }\Gamma_0.
\end{equation}

\item $\gamma<\gamma_s^*$, then it holds a  no-slip boundary condition
\begin{equation}\label{sub_bc_limit_formal}
\bar v^0=0\quad \hbox{ on }\Gamma_0.
\end{equation}
\end{itemize}
\end{theorem}
\begin{proof}  We first prove system (\ref{limit_crit_Nav}). For this, we assume the asymptotic expansion of the unknowns $(\widetilde u_\ep, \widetilde p_\ep)$ given by (\ref{formal_exp_Nav}). Then,  substituting the expansion into the problem (\ref{system_1introtilde})$_{1,2}$, we get
\begin{equation}\label{system_1introtilde_EXP}
 \begin{array}{rl}
\displaystyle
\!\!\!-\nu \ep^2\Delta_{z'}(\bar v^0+ O(\ep))- \nu\partial^2_{z_3}(\bar  v^0+ O(\ep))+
 \ep^3(v^0_3+O(\ep))\partial_{z_3}(\bar v^0+O(\ep)) +\nabla_{z'}(\widetilde  p_\ep+O(\ep))= f',\\
\\
\displaystyle
\!\!\!-\nu \ep^2\Delta_{z'}( v^0_3+ O(\ep))- \nu\partial^2_{z_3}(v^0_3+ O(\ep))+
 \ep^3(v^0_3+O(\ep))\partial_{z_3}(v^0_3+O(\ep)) +{1\over \ep}\partial_{z_3}(\widetilde  p_\ep+O(\ep))= 0,\\
\\
 \displaystyle
\!\!\! \ep^2 {\rm div}_{z'}(\bar v^0+O(\ep))+\ep \partial_{z_3}(v^0_3+\ep v^1_3+O(\ep^2))=0.
\end{array} 
\end{equation}
Collecting the terms of the same order with respect to $\ep$, we have
\begin{itemize}
\item[--] The main order terms in (\ref{system_1introtilde_EXP})$_{1,2}$ are
\begin{equation}\label{system_exp1_Nav}
 \begin{array}{rrcll}
1: & \displaystyle - \nu\partial_{z_3}^2 \bar v^0+\nabla_{z'}  p^0&=& f'& \hbox{in }\Omega,\\
\noame
{1\over \ep}: &\displaystyle \partial_{z_3}p^0&=&0& \hbox{in }\Omega.
\end{array}
\end{equation}

\item[--] The main and next order terms in (\ref{system_1introtilde_EXP})$_{3}$ are
\begin{equation}\label{system_exp2_Nav}
 \begin{array}{rrcll}
\ep : & \partial_{z_3} v^0_3&=&0 & \hbox{in }\Omega,\\
\noame
\ep^2: & {\rm div}_{x'}\bar v^0+\partial_{z_3} v^1_3&=&0& \hbox{in }\Omega,
\end{array} 
\end{equation}

\item[--] The main order and next order terms in the boundary condition $\widetilde u_\ep=0$ on  $\Gamma_1$ are 
\begin{equation}\label{bcGamma101_Nav}
 \begin{array}{rrcll}
\ep^2:& v^0&=&0& \hbox{on }\Gamma_1,\\
\noame
\ep^3:& v^1&=&0& \hbox{on }\Gamma_1.
\end{array} 
\end{equation}
\item[--] The main order and next order terms in the boundary condition $\widetilde u_{\ep,3}=0$ on $\Gamma_0$ given in  (\ref{bc1_system_tilde_0})  are 
\begin{equation}\label{bcGamma102_Nav}
 \begin{array}{rrcll}
\ep^2:& v^0_3&=&0&\hbox{on }\Gamma_0,\\
\noame
\ep^3:& v^1_3&=&0& \hbox{on }\Gamma_0.
\end{array} 
\end{equation}
\end{itemize}
From the previous equalities, we deduce
\begin{itemize}
\item[--] The limit system for $\bar v^0$  in (\ref{limit_crit_Nav})$_1$ is given in  (\ref{system_exp1_Nav})$_1$. Also, as consequence of (\ref{system_exp1_Nav})$_2$, we deduce that $p^0$ does not depend on $z_3$, i.e. $p^0=p^0(z')$.
\item[--] From (\ref{system_exp2_Nav})$_1$ and the boundary conditions (\ref{bcGamma101_Nav})$_1$ and (\ref{bcGamma102_Nav})$_1$ of $v^0_3$ on $\Gamma_0\cup \Gamma_1$, we deduce that $v_3^0 \equiv 0$.
\item[--] Integrating (\ref{system_exp2_Nav})$_2$ with respect to $z_3$ between $0$ and $h(z')$ and from the boundary conditions of $v_3^1$ on $\Gamma_1\cup \Gamma_0$ given respectively in (\ref{bcGamma101_Nav})$_2$ and (\ref{bcGamma102_Nav})$_2$,   we deduce the incompressibility condition (\ref{limit_crit_Nav})$_2$.
\end{itemize}
To finish the proof, it remains to deduce the boundary conditions of $\bar v^0$ on $\Gamma_0$.  Multiplying by $\ep^{-1}$ the boundary condition (\ref{bc1_system_tilde_0}) 
and using the expansion (\ref{formal_exp_Nav}), we have
$$
-{1\over \ep^2}\nu\partial_{z_3}(\ep^2 \bar v^0+ \ep^3 \bar v^1+\cdots)=\ep^{\gamma-1}|K(\ep^2 \bar v^0+ \ep^3 \bar v^1+\cdots)|^{s-2}K^2(\ep^2 \bar v^0+ \ep^3 \bar v^1+\cdots),
$$
 which, after simplification, implies 
 \begin{equation}\label{expansion_partial_condition}
-\nu\partial_{z_3}(\bar v^0+ \ep \bar v^1+\cdots)=\ep^{\gamma+1}|K(\ep^2 \bar v^0+ \ep^3 \bar v^1+\cdots)|^{s-2}K^2( \bar v^0+ \ep \bar v^1+\cdots).
\end{equation}
 Taking into account that  
$$|K(\ep^2 \bar v^0+ \ep^3 \bar v^1+\cdots)|^{s-2}=\ep^{2s-4}|K(\bar v^0+ \ep \bar v^1+\cdots)|^{s-2}=\ep^{2s-4}|K\bar v^0+ \mathcal{O}(\ep)|^{s-2},$$
then, we have that the main order term in (\ref{expansion_partial_condition}) is
 \begin{equation}\label{expansion_partial_condition2}
-\nu\partial_{z_3}\bar v^0=\ep^{\gamma-\gamma_s^*}|K \bar v^0|^{s-2}K^2 \bar v^0.
\end{equation}
Thus, we have:
\begin{itemize}
\item If $\gamma=\gamma_s^*$, we deduce the boundary condition (\ref{crit_bc_limit_formal}).

\item If $\gamma>\gamma_s^*$, for example, we can consider in (\ref{expansion_partial_condition2})  the following value of $\gamma$,
$$\gamma=\gamma_s^*+k,\quad k\in \mathbb{Z}_+.$$ 
Then, $\gamma-\gamma_s^*>0$ and so, the main order term in (\ref{expansion_partial_condition2}) is
$$-\nu\partial_{z_3}\bar v^0=0,$$ 
which is the boundary condition (\ref{super_bc_limit_formal}).

\item If $\gamma<\gamma_s^*$, for example, we can consider in (\ref{expansion_partial_condition2})  the following value of $\gamma$,
$$\gamma=\gamma_s^*-k,\quad k\in \mathbb{Z}_+,$$ 
Then, $\gamma-\gamma_s^*<0$ and so, the main order term in (\ref{expansion_partial_condition2}) is
$$|K \bar v^0|^{s-2}K^2 \bar v^0=0,$$
which implies the boundary condition (\ref{sub_bc_limit_formal}).
\end{itemize}

\end{proof}

\section{Mathematical justification} \label{sec:Rigorous}
In this section, we derive {\it a priori} estimates for $(\widetilde u_\ep, \widetilde p_\ep)$, which will justify the formal asymptotic expansion (\ref{formal_exp_Nav}), and we give the corresponding compactness results leading to the limit problem, which is the main result of this paper.
\subsection{A priori  estimates}
We start  deducing a priori estimates for the rescaled functions $\widetilde u_\ep $.  Let us first recall some important inequalities in a thin domain with thickness $\ep$, see for instance \cite{AnguianoSGNonlinear, AnguianoSGmag, MikTap}.
\begin{lemma}[Poincar\'e's inequality]\label{Poincare_lemma} For all $\varphi\in \mathbb{V}(\Omega^\ep)$,  the following inequality holds
\begin{equation}\label{Poincare}
\|\varphi\|_{L^2(\Omega^\varepsilon)^3}\leq C\varepsilon\|D \varphi\|_{L^2(\Omega^\varepsilon)^{3\times 3}}.
\end{equation}
Moreover, from the change of variables (\ref{dilatation}),   there hold the following rescaled estimate
\begin{equation}\label{Poincare2}
\|\widetilde \varphi\|_{L^2( \Omega)^3}\leq C \varepsilon\|D_{ \ep} \widetilde \varphi\|_{L^2( \Omega)^{3\times 3}}.
\end{equation}
\end{lemma}
 Let us introduce the version of the Ladyzhenskaya's inequality, which enable us to estimate the inertial term in the variational formulation.
\begin{lemma}[Ladyzhenskaya's inequality]\label{Lady_lemma} For all $\varphi\in H^1_0(\Omega^\ep)^3$,  the following inequality hold
\begin{equation}\label{Lady}
 \|\varphi\|_{L^4(\Omega^\varepsilon)^3}\leq 
 C\varepsilon^{1\over 4}\|D \varphi\|_{L^2(\Omega^\varepsilon)^{3\times 3}}.
\end{equation}
Moreover, from the change of variables (\ref{dilatation}),   there hold the following rescaled estimate
\begin{equation}\label{Lady2}
\|\widetilde \varphi\|_{L^4(\Omega)^3}\leq C\varepsilon^{1\over 2}\|D_{\ep} \widetilde \varphi\|_{L^2( \Omega)^{3\times 3}}.
\end{equation}
\end{lemma}
\begin{proof}
Using the interpolation inequality, the embedding $H^1(\Omega^\ep)\hookrightarrow L^6(\Omega^\ep)$,  and the estimate (\ref{Poincare}), we obtain at once
$$
\|\varphi\|_{L^4(\Omega^\ep)^3}\leq \|\varphi\|_{L^2(\Omega^\ep)^3}^{1\over 4}\|\varphi\|_{L^6(\Omega^\ep)^3}^{3\over 4}\leq C\ep^{1\over 4}\|D\varphi\|_{L^2(\Omega^\ep)^{3\times 3}}.$$
Estimate (\ref{Lady2}) follows from the change of variables  (\ref{dilatation}) by taking into account that 
$$\|\varphi\|_{L^4(\Omega^\varepsilon)^3}=\ep^{1\over 4}\|\widetilde \varphi\|_{L^4(\Omega)^3},\quad \|D \varphi\|_{L^2(\Omega^\varepsilon)^{3\times 3}}=\ep^{1\over 2}\|D_{\ep} \widetilde \varphi\|_{L^2( \Omega)^{3\times 3}}.$$
\end{proof}

Now, we prove a trace estimate for the rescaled velocity.
\begin{lemma}[Trace estimates] For all $\widetilde\varphi \in \mathbb{V}(\Omega)$, the following inequalities hold\begin{equation}\label{trace_estimate}
\|\widetilde\varphi\|_{L^2(\Gamma_0)^3}\leq C\|D_z\widetilde \varphi\|_{L^2(\Omega)^{3\times 3}},\quad \|\widetilde\varphi\|_{L^2(\Gamma_0)^3}\leq  C\ep\|D_\ep \widetilde \varphi\|_{L^2(\Omega)^{3\times 3}}.
\end{equation}
\end{lemma}
\begin{proof}
Thank to $\widetilde \varphi(z',h(z'))=0$ in $\omega$, we have that
$$\begin{array}{rl}
\displaystyle \int_{\Gamma_0}|\widetilde \varphi|^2d\sigma= &\displaystyle \int_{\omega}|\widetilde \varphi(z',0)|^2\,dz'= \int_{\omega}\left|\int_0^{h(z_1)}\partial_{z_3}\widetilde \varphi(z)\,dz_3\right|^2\,dz'\leq C  \int_{\Omega}|\partial_{z_3}\widetilde \varphi(z)|^2\,dz',
\end{array}$$
that is, 
$$\|u'_\ep\|_{L^2(\Gamma_0)^3}\leq C\|\partial_{z_3}\widetilde \varphi\|_{L^2(\Omega)^{3\times 3}},$$
which implies (\ref{trace_estimate})$_1$ and  (\ref{trace_estimate})$_2$.

\end{proof}

\noindent Next, we give the {\it a priori} estimates of velocity.

\begin{lemma}\label{lem:estimvelmic} Assume $1<s<2$ and let $\gamma\in \mathbb{R}$.  Let $u_\ep$ be a weak solution of (\ref{system_1introad})--(\ref{bc0_system}). 
Then, there exists a positive constant $C$, independent of $\ep$, such that we have the following estimates
\begin{equation}\label{velocityestim_power}
\| u_\ep \|_{L^2(\Omega^\ep  )^3}\leq C\ep^{5\over 2},\quad \|D  u_\ep \|_{L^2(\Omega^\ep  )^{3\times 3}}\leq C\ep^{3\over 2} ,
\end{equation}
\begin{equation}\label{velocityestim_boundary}
\| K u_\ep' \|_{L^s(\Gamma_0)^2}\leq C\ep^{3-\gamma\over s}.
\end{equation}
Moreover, after the change of variables (\ref{dilatation}), we have the following estimates
\begin{equation}\label{velocityestimtilde_power}
\| \widetilde u_\ep \|_{L^2(\Omega)^3}\leq C\ep^{2},\quad \|D_\ep \widetilde u_\ep \|_{L^2(\Omega )^{3\times 3}}\leq C\ep,
\end{equation}
\begin{equation}\label{velocityestim_boundary_tilde}
\| K\widetilde u_\ep' \|_{L^s(\Gamma_0)^2}\leq C\ep^{3-\gamma\over s}.
\end{equation}
 
\end{lemma}
\begin{proof}  
From (\ref{form_var_vel}) with $\varphi=  u_\ep $, we get
\begin{equation}\label{estim_vel}
\begin{array}{c}
\displaystyle
\nu \|D u_\ep \|_{L^2(\Omega^\ep)^{3\times 3}}^2 +\ep^{\gamma} \|Ku_\ep'\|_{L^s(\Gamma_0)^2}^s =\int_{\Omega^\ep }f'\cdot  u_\ep'\,dx,
\end{array}
\end{equation}
because $\int_{\Omega^\ep}(u_\ep \cdot \nabla)  u_\ep \,  u_\ep \,dx=0$.  From the Cauchy-Schwarz  and Poincar\'e's inequality (\ref{Poincare}), it holds
$$\int_{\Omega^\ep}f'\cdot \widetilde u_\ep'\,dx\leq C\ep^{1\over 2}\|u^\ep\|_{L^2(\Omega^\ep)^3}\leq C\ep^{3\over 2}\|D u^\ep\|_{L^2(\Omega^\ep)^{3\times 3}}.$$
This, together with (\ref{estim_vel}), implies the second estimate in (\ref{velocityestim_power}). Next, by using again Poincar\'e's inequality (\ref{Poincare}), we deduce the first one.   On the other hand, from (\ref{estim_vel}) and (\ref{velocityestim_power})$_2$, then 
$$\ep^{\gamma} \|Ku_\ep'\|_{L^s(\Gamma_0)^2}^s\leq C\ep^3,$$
and so,  estimate (\ref{velocityestim_boundary}) holds.\\

\noindent Finally, by applying the dilatation (\ref{dilatation}) to previous estimates (\ref{velocityestim_power}) and (\ref{velocityestim_boundary}), and taking into account 
taking into account 
$$\ep^{1\over 2}\| \widetilde u_\ep \|_{L^2(\Omega)^3}=\| u_\ep \|_{L^2(\Omega^\ep)^3},\quad \ep^{1\over 2}\|D_\ep \widetilde u_\ep \|_{L^2(\Omega )^{3\times 3}}=\|D u_\ep \|_{L^2(\Omega^\ep )^{3\times 3}},\quad \| K\widetilde u_\ep' \|_{L^s(\Gamma_0)^2}=\|Ku_\ep' \|_{L^s(\Gamma_0)^2},$$
we deduce estimates (\ref{velocityestimtilde_power}) and (\ref{velocityestim_boundary_tilde}), respectively.\\
\end{proof}

Now, we will derive the estimate of the  pressure by using the {\it a priori} estimates for velocity given in Lemma \ref{lem:estimvelmic}.
\begin{lemma}\label{lem:estimpress}   Let  $\widetilde p_\ep$ be a weak solution of problem  (\ref{system_1introad})--(\ref{bc0_system}). Then, there exists a positive constant $C$, independent of $\ep$, such that  we have the following estimates 
\begin{equation}\label{estim_p1}
\begin{array}{c}
\displaystyle
\|\widetilde p_\varepsilon\|_{L^{2}_0(\Omega)}\leq C,\quad \|\nabla_\ep \widetilde p_\ep \|_{H^{-1}(\Omega)^3}\leq C.
\end{array}
\end{equation}
\end{lemma}
\begin{proof} Thank to the classical Ne${\rm \check{c}}$kas estimate  (see \cite{Boyer})
\begin{equation}\label{Neckas}
\|\widetilde p_\ep \|_{L^{2}_0(\Omega)}\leq C\|\nabla_z \widetilde p_\ep \|_{H^{-1}(\Omega)^3},
\end{equation}
and taking into account that $\|\nabla_z \widetilde p_\ep \|_{H^{-1}(\Omega)^3}\leq C\|\nabla_\ep \widetilde p_\ep \|_{H^{-1}(\Omega)^3}$, then we just need to obtain the estimate for $\nabla_\ep \widetilde p_\ep $ to derive the estimates for $\widetilde p_\ep $ given in (\ref{estim_p1}) with the restriction imposed on the possible values of $\gamma$.\\

\noindent To do this, we consider $\widetilde \varphi\in H^1_0(\Omega)^3$, and taking into account the variational formulation  (\ref{form_var_vel_tilde}), we get
\begin{equation}\label{equality_duality_0}
\begin{array}{rl}
\displaystyle
\left\langle \nabla_\ep  \widetilde p_\varepsilon, \widetilde \varphi\right\rangle_{H^{-1}(\Omega)^3,H^1_0(\Omega)^3}=&\displaystyle
-  \nu\int_{\Omega}D_\ep\widetilde u_\ep :D_\ep\widetilde \varphi\,dz -\int_{\Omega}(\widetilde u_\ep\cdot \nabla_\ep)\widetilde u_\ep\cdot \widetilde\varphi\,dz'+\int_{\Omega}f'\cdot   \tilde \varphi'\,dz.
\end{array}\end{equation}

\noindent We estimate  the terms on the right-hand side of  (\ref{equality_duality_0}):
\begin{itemize}
\item[--] First term. Applying  Cauchy-Schwarz's inequality and estimate  (\ref{velocityestimtilde_power}), we get
$$\begin{array}{rcl}
\displaystyle\left|\nu\int_{\Omega}D_\ep\widetilde u_\ep :D_\ep\widetilde \varphi\,dz\right| &\leq &\displaystyle
C\|D_\ep \widetilde  u_\ep \|_{L^2( \Omega)^{3\times 3}}\|D_\ep \widetilde \varphi\|_{L^2(\Omega)^{3\times 3}}\\
\noame
&\leq &   C\ep^{-1}\|D_\ep \widetilde  u_\ep \|_{L^2( \Omega)^{3\times 3}}\|\widetilde \varphi\|_{H^1_0(\Omega)^{3}}\\
\noame
&\leq& C\|\widetilde \varphi\|_{H^1_0(\Omega)^{3}}.
\end{array}$$

\item[--] Second term. Applying H${\rm \ddot{o}}$lder's inequality and using Ladyzhenskaya's inequality  (\ref{Lady2})  and  (\ref{velocityestimtilde_power}), we have  
$$\begin{array}{rcl}
\displaystyle
 \left|\int_{\Omega}(\widetilde u_\ep \cdot \nabla_\ep)\widetilde u_\ep \cdot \widetilde \varphi\,dz\right| &\leq &\displaystyle  
 C \|\widetilde  u_\ep \|_{L^4(\Omega)^{3}}\|D_\ep  \widetilde  u_\ep \|_{L^2( \Omega)^{3\times 3}}\|\widetilde  \varphi\|_{L^4(\Omega)^{3}}\\
 \noame 
 &\leq &   C \ep \|D_\ep  \widetilde  u_\ep \|^2_{L^2( \Omega)^{3\times 3}} \|D_\ep \widetilde  \varphi\|_{L^2( \Omega)^{3\times 3}}\\
 \noame 
 &\leq &   C \ep^2 \| \widetilde  \varphi\|_{H^1_0( \Omega)^{3}}.
\end{array}$$

\item[--]  Third term. Applying Cauchy-Schwarz's inequality,  $f'\in L^2(\omega)^2$ and Poincar\'e's inequality (\ref{Poincare2}), we get
$$\begin{array}{rcl}
\displaystyle \left|\int_{\Omega}f'\cdot  \widetilde \varphi'\,dz\right|
& \leq &\displaystyle 
C\|\widetilde \varphi\|_{L^2(\Omega)^3}\leq C \varepsilon\|D_{ \ep} \widetilde \varphi\|_{L^2( \Omega)^{3\times 3}}\leq  C\|\widetilde \varphi\|_{H^1_0(\Omega)^3}.
\end{array}$$
\end{itemize}

\noindent Finally, taking into account all the previous estimates, we get the desired estimate
$$\|\nabla_\ep \widetilde p_\ep\|_{H^{-1}(\Omega)^3}\leq 
C,$$
which finishes the proof.
 
\end{proof}

\subsection{Compactness results}\label{subsec:main1}

 First, we give some compactness results about the behavior of the rescaled sequences $\widetilde u_\ep$ and $\widetilde p_\ep $ satisfying respectively the a priori estimates given in Lemmas \ref{lem:estimvelmic} and \ref{lem:estimpress}. Next, we deduce the limit system depending on the case and the value of $\gamma$, which is the main result of this paper (see Theorem \ref{thm1_critical_rig}).
 
\begin{lemma} \label{lemmaconv_crit}  Assume $1<s<2$, $\gamma\in \mathbb{R}$ and $\gamma_s^*$ given by (\ref{gammastar}).  There exist:
\begin{itemize}
 \item A subsequence, still denoted by $(\widetilde u_\ep, \widetilde p_\ep)$, of the sequence of  solutions $(\widetilde u_\ep, \widetilde p_\ep)$ of system  (\ref{system_1introad})--(\ref{bc0_system}).

\item   $\widetilde u \in V_{z_3}$ and $\widetilde u_3\equiv 0$,  such that 
\begin{equation}\label{velconv_power}
\ep^{-2}\widetilde u_\ep \rightharpoonup \widetilde u\quad\hbox{in}\quad V_{z_3},
\end{equation}
\begin{equation}\label{velconvdiv_power}
{\rm div}_{z'}\left(\int_0^{h(z')}\widetilde u'(z)\,dz_3\right)=0\quad\hbox{in }  \omega,\quad \left(\int_0^{h(z')}\widetilde u'(z)\,dz_3\right)\cdot {\rm n} =0\quad\hbox{on } \partial\omega,
\end{equation}
where $V_{z_3}$ is the Hilbert space (see \cite{Bayada}) 
$$V_{z_3}=\{\varphi\in L^2(\Omega)^3\, :\, \partial_{z_3}\varphi\in L^2(\Omega)^3\},$$
with the norm
$$\|\varphi\|_{V_{z_3}}^2=\|\varphi\|_{L^2(\Omega)^3}^2+\|\partial_{z_3}\varphi\|^2_{L^2(\Omega)^3}.$$
Moreover,   the following values of the velocity $\widetilde u'$ on the boundaries hold:
\begin{itemize}
\item If $\gamma\geq \gamma_s^*$, then $\widetilde u'$ satisfies no-slip condition on the top boundary, i.e.  \begin{equation}\label{noslipG1_power}\widetilde u'=0\quad \hbox{on}\quad \Gamma_1.
\end{equation}

\item If $\gamma< \gamma_s^*$, then $\widetilde u'$ satisfies no-slip condition on the top and bottom boundaries, i.e.  \begin{equation}\label{noslipG10_power}\widetilde u'=0\quad \hbox{on}\quad \Gamma_0\cup \Gamma_1.
\end{equation}
\end{itemize}

\item $\tilde p \in L^{2}_0(\Omega)$, independent of $z_3$,  such that 
\begin{equation}\label{pressconv_power}
\widetilde p_\ep \rightharpoonup \widetilde p\quad\hbox{in}\quad L^{2}(\Omega).
\end{equation}
\end{itemize}

\end{lemma}
\begin{proof} The convergences given in this lemma are a direct consequence of the a priori estimates given in Lemmas \ref{lem:estimvelmic} and \ref{lem:estimpress}. We will only give some remarks:
\begin{itemize}
\item[--] {\it Velocity}. Estimates given in   (\ref{velocityestimtilde_power}) imply the existence of a function $\widetilde u\in V_{z_3}$ such that convergence   (\ref{velconv_power}) holds. The continuity of the trace applications from the space of $\widetilde u$ such that $\|\widetilde u\|_{L^2}$ and $\|\partial_{z_3}\widetilde u\|_{L^2}$ are bounded to   $L^2(\Gamma_1)$ implies that $\widetilde u=0$ on $\Gamma_1$. Also, since the continuity is also bounded to $L^2(\Gamma_0)$, from $u_{\ep,3}=0$ on $\Gamma_0$ it holds that $\widetilde u_3=0$ on $\Gamma_0$. \\

From the variational formulation (\ref{formdiv_tilde})  for  $\psi\in C^1(\Omega)$, after multiplication by $\ep^{-1}$, integrating by parts and taking into account that $\widetilde u_{\ep}\cdot {\rm n}=0$ on $\partial\Omega$, we have 
$$
\int_{\Omega}\ep^{-1} \widetilde u_\ep\cdot\nabla_\ep \widetilde \psi\,dz=0,\ \hbox{i.e. }\quad \int_{\Omega}\left(\ep^{-1} \widetilde u_\ep'\cdot\nabla_{z'} \widetilde \psi+\ep^{-2}\widetilde u_{\ep, 3}\partial_{z_3}\widetilde\psi\,\right)dz=0.
$$
Taking into account estimate    (\ref{velocityestimtilde_power}),  it can be written as follows
$$
\int_{\Omega}\ep^{-2} \widetilde u_{\ep,3}\partial_{z_3}\psi\,dz+O_\ep=0,
$$
and thus, from convergence   (\ref{velconv_power}), we deduce
$$
\int_{\Omega}  \widetilde u_3\,\partial_{z_3}\psi\,dz=0,
$$
which implies  $\partial_{z_3} \widetilde u_3=0$,  i.e. $\widetilde u_3$ is independent of $z_3$.  This combined with the boundary condition  $\widetilde u_3=0$ on $ \Gamma_0\cup \Gamma_1$, implies $\widetilde u_3\equiv 0$. 
\\

Now, from the variational formulation (\ref{formdiv_tilde})    for   $\psi\in C^\infty(\omega)$, we have 
$$
\int_{\Omega}\widetilde u_\ep'\cdot\nabla_{z'}\widetilde \psi'\,dz=0,
$$
and multiplying by $\ep^{-2}$ and using convergence   (\ref{velconv_power}), we deduce
$$
\int_{\Omega}\widetilde u'\cdot\nabla_{z'}\widetilde \psi'\,dz=\int_{\omega}{\rm div}_{z'}\left(\int_0^{h(z')}\widetilde u'\,dz_3\right) \widetilde \psi'\,dz'=0.
$$
This implies the divergence condition (\ref{velconvdiv_power}).\\

It remains to prove  that $\widetilde u'=0$ on $\Gamma_0$ in the case $\gamma<\gamma_s^*$. To do this, from the continuous embedding of $V_{z_3}$ into $L^2(\Gamma_0)$, we deduce from  (\ref{velocityestimtilde_power})  that  
$$\|K\widetilde u'_\ep\|_{L^2(\Gamma_0)^2}\leq \|K\|_{L^\infty(\Gamma_0)^{2\times 2}}\|\widetilde u'_\ep\|_{L^2(\Gamma_0)^2}\leq C \ep^2,$$
  and   since $1<s<2$, it holds
$$  \|K\widetilde u'_\ep\|_{L^s(\Gamma_0)^2}\leq  C \|K\widetilde u'_\ep\|_{L^2(\Gamma_0)^2}\leq C\ep^2.$$
This implies
$$K(\ep^{-2}\widetilde u'_\ep)\rightharpoonup K\widetilde u' \quad\hbox{in }L^s(\Gamma_0)^2.$$
On the other hand, since we have
$$K(\ep^{-2}\widetilde u'_\ep)=\ep^{{3-\gamma\over s}-2}(K\ep^{\gamma-3\over s}\widetilde u'_\ep)=\ep^{{\gamma_s^*-\gamma\over s}} (K\ep^{\gamma-3\over s}\widetilde u'_\ep),$$
then, from estimate (\ref{velocityestim_boundary_tilde}) and $\ep^{{\gamma_s^*-\gamma\over s}} \to 0$, we get that 
$$K(\ep^{-2}\widetilde u'_\ep)\rightharpoonup 0\quad\hbox{in }L^s(\Gamma_0)^2.$$
Thus, from the uniqueness of the limit and $K$ is positive definite, we deduce $\widetilde u'=0$ on $\Gamma_0$.

\item[--] {\it Pressure}. Estimate (\ref{estim_p1})$_1$ implies, up to a subsequence, the existence of $\widetilde p\in L^2(\Omega)$  such that convergence  (\ref{pressconv_power}) holds. Also, from (\ref{estim_p1})$_2$, by noting that $\ep^{-1}\partial_{z_3}\widetilde p_\ep$ also converges weakly in $H^{-1}(\Omega)$, we deduce $\partial_{z_3}\widetilde p=0$. Then, $\widetilde p$  is independent of $z_3$. To finish, it remains to prove that $p \in L^2_0(\Omega)$. Passing to the limit when $\ep$ tends to zero in 
$$\int_{\Omega}\widetilde p_\ep\,dz=0,$$
we respectively deduce
$$\int_{\Omega}\widetilde p(z')\,dz=\int_{\omega}h(z')\widetilde p(z')\,dz'=0,$$
and so that $\widetilde p$  has null mean value in $\Omega$. This ends the proof.
\end{itemize}
\end{proof}
Next, we prove the main result of this paper.
\begin{theorem}[Main result] \label{thm1_critical_rig}
The limit pair of functions $(\widetilde u, \widetilde p)\in V_{z_3}\times L^2_0(\Omega)$, with $\widetilde u_3\equiv 0$  and $\widetilde p=\widetilde p(z')$, given in Lemma \ref{lemmaconv_crit} satisfies the following system in each case:
\begin{equation}\label{limit_crit_rig}
\left\{\begin{array}{rl}
\displaystyle - \nu\partial_{z_3}^2\widetilde u' (z)=  f'(z')-\nabla_{z'}\widetilde p(z') & \hbox{in } \Omega,\\
\noame
\displaystyle{\rm div}_{z'}\left(\int_0^{h(z')}\widetilde u'(z)\,dz_3\right)=0& \hbox{in } \omega,\\
\noame
\displaystyle \left(\int_0^{h(z')}\widetilde u'(z)\,dz_3\right)\cdot {\rm n} =0& \hbox{on }\partial\omega,\\
\noame
\displaystyle \widetilde u'=0&\hbox{on }\Gamma_1,
\end{array} \right.
\end{equation}
with $\nu>0$. Moreover, $\widetilde u'$ satisfies the following boundary condition on the bottom $\Gamma_0$ depending on the value of $\gamma$:
\begin{itemize}
\item $\gamma=\gamma_s^*$, then it holds a power law slip boundary condition
\begin{equation}\label{crit_bc_limit}
 -\nu\partial_{z_3}\widetilde u'=-|K\widetilde u'|^{s-2}K^2\widetilde u'\quad \hbox{ on }\Gamma_0,
\end{equation}
where   $1<s<2$ and the anisotropic tensor $K\in\mathbb{R}^{2\times 2}$  is  uniformly positive definite, symmetric and bounded.
\item $\gamma>\gamma_s^*$, then it holds a perfect slip boundary condition
\begin{equation}\label{super_bc_limit}
 -\nu\partial_{z_3}\widetilde u'=0\quad \hbox{ on }\Gamma_0.
\end{equation}

\item $\gamma<\gamma_s^*$, then it holds a  no-slip boundary condition
\begin{equation}\label{sub_bc_limit}
\widetilde u'=0\quad \hbox{ on }\Gamma_0.
\end{equation}

\end{itemize}

\end{theorem}
\begin{remark}
By uniqueness of solutions of problems given in  Theorem \ref{thm1_critical_rig}, we observe that the pair of functions $(\widetilde u,\widetilde p)$ are the same as those functions $(v^0, p^0)$ obtained in Theorem \ref{thm:formal_NS} by formal arguments.
\end{remark}
\begin{proof}[Proof of Theorem \ref{thm1_critical_rig}] We will divide the proof in three steps. \\

{\it Step 1}. Let us first consider the case $\gamma>\gamma_s^*$. From Lemma \ref{lemmaconv_crit}, to prove (\ref{limit_crit_rig}), it just remains to prove the equations (\ref{limit_crit_rig})$_1$.  To do this,  we  consider  $\widetilde\varphi\in C^1_c(\omega\times (0,h(z')))^3$ such that $\widetilde\varphi_3\equiv 0$ and $\widetilde \varphi=0$ on $\Gamma_1$. Thanks to $\widetilde\varphi$ equaling zero for $z'$ outside a compact subset of $\omega$, then
 $$\widetilde\varphi=0\ \hbox{on}\ \partial\Omega\setminus\Gamma_0\quad \hbox{and}\quad\widetilde \varphi_3=0 \ \hbox{on}\ \Gamma_0,$$
 and so, we can take it  as test function in (\ref{form_var_vel_tilde}), which is given by
\begin{equation}\label{form_var_veltildeproof1}
\begin{array}{l}
\displaystyle
 \nu \int_{\Omega}D_\ep \widetilde u_\ep :D_\ep\widetilde \varphi\,dz+\int_{\Omega}(\widetilde u_\ep\cdot \nabla_\ep)\widetilde u_\ep\cdot \widetilde\varphi\,dz+\ep^{\gamma-1} \int_{\Gamma_0}|K\widetilde u_\ep'|^{s-2}K \widetilde u_\ep'\cdot  K\widetilde \varphi'\,d\sigma
 \\
 \noame
 \displaystyle -\int_{\Omega}\widetilde p_\ep\,{\rm div}_\ep(\widetilde \varphi)\,dz=\int_{\Omega}f'\cdot \widetilde \varphi'\,dz.
 \end{array}
\end{equation}

\noindent Let us now pass to the limit when $\ep$ tends to zero in every terms of (\ref{form_var_veltildeproof1}):
\begin{itemize}
\item[--] First term on the left-hand side. Taking into account convergence (\ref{velconv_power}), we get
$$  \nu \int_{\Omega}D_\ep \widetilde u_\ep :D_\ep\widetilde \varphi\,dz= \nu \int_{\Omega}\ep^{-2}\partial_{z_3}\widetilde u_{\ep}'\cdot \partial_{z_3}\widetilde \varphi'\,dz+O_\ep= \nu\int_{\Omega} \partial_{z_3}\widetilde u'\cdot \partial_{z_2}\widetilde \varphi'\,dz+O_\ep.$$
\item[--] Second term on the left-hand side. Taking into account the regularity of $\varphi'$, applying Cauchy-Schwarz's inequality and estimates (\ref{velocityestimtilde_power}), we get
 \begin{equation}\label{inertial_variational}\left|\int_{\Omega}(\widetilde u_\ep \cdot \nabla_\ep)\widetilde u_\ep \widetilde \varphi\,dz\right|\leq \|\widetilde u_\ep \|_{L^2(\Omega)^3}\|D_\ep \widetilde u_\ep \|_{L^2(\Omega)^{3\times 3}}\|\widetilde \varphi\|_{L^\infty(\Omega)^3}\leq C\ep^{3},
 \end{equation}
which implies
$$\int_{\Omega}(\widetilde u_\ep \cdot \nabla_\ep)\widetilde u_\ep\cdot  \widetilde \varphi\,dz\to 0.$$
\item[--] Third term on the left-hand side. We observe that since $\Gamma_0$ is flat, then the surface measure associated to $\Gamma_0$ given by $d\sigma=dz'$. From H${\rm \ddot{o}}$lder's inequality, $K$ is bounded, the Sobolev embedding $L^2\hookrightarrow L^s$, the trace estimate (\ref{trace_estimate})$_1$ applied to $\widetilde \varphi'$, the trace estimate (\ref{trace_estimate})$_2$ applied to $\widetilde u_\ep'$, and the estimate (\ref{velocityestimtilde_power}) for $D_\ep\widetilde u_\ep'$, we get
$$
\begin{array}{rl}
\displaystyle \left|\ep^{\gamma-1} \int_{\Gamma_0}|K\widetilde u_\ep'|^{s-2} K\widetilde u_\ep'\cdot  K\widetilde \varphi'\,d\sigma\right|  \leq &\displaystyle \ep^{\gamma-1}\|K\widetilde u_\ep'\|_{L^s(\Gamma_0)^2}^{s-1}\|K\widetilde \varphi'\|_{L^s(\Gamma_0)^2}\\
\noame
\leq &\displaystyle 
\ep^{\gamma-1}\|K\|_{L^\infty(\Gamma_0)^{2\times 2}}^s\|\widetilde u_\ep'\|_{L^s(\Gamma_0)^2}^{s-1}\|\widetilde \varphi'\|_{L^s(\Gamma_0)^2}
\\
\noame
\leq &\displaystyle 
C\ep^{\gamma-1}\|\widetilde u_\ep'\|_{L^2(\Gamma_0)^2}^{s-1}\|\widetilde \varphi'\|_{L^2(\Gamma_0)^2}
\\
\noame
\leq&\displaystyle
  C\ep^{\gamma-1}\ep^{s-1}\|D_\ep \widetilde u_\ep'\|_{L^2(\Gamma_0)^{3\times 2}}^{s-1} \|D_z\widetilde \varphi'\|_{L^2(\Omega)^{3\times 2}}\\
\noame
\leq&\displaystyle
 C\ep^{\gamma-\gamma_s^*}\|D_z\widetilde \varphi'\|_{L^2(\Gamma_0)^{3\times 2}}\\
 \noame
 \leq &\displaystyle C\ep^{\gamma-\gamma_s^*},
 \end{array}
$$
which tends to zero because $\gamma>\gamma_s^*$. Then, we get 
$$\ep^{\gamma-1} \int_{\Gamma_0}|K\widetilde u_\ep'|^{s-2} K\widetilde u_\ep'\cdot K \widetilde \varphi'\,d\sigma\to 0.$$
\item[--] Fourth term on the left-hand side of (\ref{form_var_veltildeproof1}). Taking into account that $\widetilde\varphi_3\equiv 0$ and convergence (\ref{pressconv_power}), we get
$$\int_{\Omega}\widetilde p_\ep \,{\rm div}_\ep(\widetilde \varphi)\,dz=\int_{\Omega}\widetilde p_\ep \,{\rm div}_{x'}(\widetilde \varphi')\,dz=\int_{\Omega}\widetilde p\,{\rm div}_{x'}(\widetilde \varphi')\,dz+O_\ep.
$$
\end{itemize}
Finally, from the above convergences when $\ep\to 0$, we derive the following limit system 
\begin{equation}\label{form_limit_Step1}\begin{array}{l}\displaystyle 
\nu\int_{\Omega} \partial_{z_3}\widetilde u'\cdot \partial_{z_3}\widetilde \varphi'\,dz-\int_{\Omega}\widetilde p\,{\rm div}_{z'}(\widetilde \varphi')\,dz=
\int_{\Omega}f'\cdot \widetilde\varphi'\,dz,
\end{array}
\end{equation}
for every $\widetilde\varphi'\in C^1_c(\omega\times (0,h(z')))^2$ with $\varphi'=0$ on $\Gamma_1$. By density, this equality holds true for every $\widetilde\varphi'\in H^1(0,h(z');L^2(\omega)^2)$ such that $\widetilde \varphi'=0$ on $\Gamma_1$. We observe that problem (\ref{form_limit_Step1}) has a unique solution  $(\widetilde u', \widetilde p)$ and the problem is equivalent to (\ref{limit_crit_rig})$_1$ with boundary condition (\ref{super_bc_limit}). Uniqueness of solution for (\ref{limit_form_var_K}) implies that limit does not depend on the subsequence.
\\

{\it Step 2}. Next, we consider the case $\gamma<\gamma_s^*$. According to Lemma \ref{lemmaconv_crit}, we proceed similarly to the Step 1, but here we also consider $\widetilde\varphi=0$ on $\Gamma_0$, i.e. we consider a test function $\widetilde\varphi\in C^1_c(\Omega)^3$ with $\widetilde\varphi_3\equiv 0$. This means that there is no boundary term  in the variational formulation (\ref{form_var_veltildeproof1}). Thus, proceeding as  Step 1, we deduce the limit variational formulation (\ref{form_limit_Step1}), which holds for every $\widetilde\varphi'\in H^1_0(\Omega)^2$. This problem has a unique solution and is equivalent to problem  (\ref{limit_crit_rig})$_1$ with boundary condition (\ref{sub_bc_limit}). Uniqueness of solution for (\ref{limit_form_var_K}) implies that limit does not depend on the subsequence.
 \\

{\it Step 3}. Finally, we consider the case $\gamma=\gamma_s^*$. Due to the nonlinear boundary term in (\ref{form_var_veltildeproof1}), we need to use monotonicity arguments to pass to the limit. For this, to simplify the notation,   we define  the application $\widetilde v\mapsto\mathcal{A}_\ep(\widetilde v) $  as follows
$$(\mathcal{A}_\ep(\widetilde v), \widetilde w)=\nu\int_{\Omega}D_\ep \widetilde v:D_\ep \widetilde w\,dz+\ep^{\gamma-1}\int_{\Gamma_0}|K \widetilde v'|^{s-2}K\widetilde v'\cdot K\widetilde w'\,d\sigma,$$
for all $\widetilde v, \widetilde w\in H^1(\Omega)^3$ such that $\widetilde \varphi=\widetilde w=0$ on $\partial\Omega\setminus \Gamma_0$ and $\widetilde \varphi_3\equiv\widetilde w_3 \equiv 0$ on $\Gamma_0$.  From \cite[Lemma 2.3]{Djoko}, for every $\ep>0$, the mapping $\mathcal{A}_\ep$ is strictly monotone, i.e. 
\begin{equation}\label{Amonotone}
(\mathcal{A}_\ep(v)-\mathcal{A}_\ep(w), v-w)\geq 0.
\end{equation}
According to Lemma \ref{lemmaconv_crit} and similar to Step 1, we consider $\widetilde\varphi\in C^1_c(\omega\times (0,h(z')))^3$ with $\widetilde\varphi_3\equiv 0$ and $\varphi=0$ on $\Gamma_1$, and we choose $\widetilde  v_\ep$ defined by
$$\widetilde v_\ep=\widetilde \varphi- \widetilde u_\ep,$$
as test function in (\ref{form_var_veltildeproof1}). So we get
$$
\begin{array}{l}
\displaystyle
(\mathcal{A}_\ep(\widetilde u_\ep),\widetilde v_\ep) -\int_{\Omega}\widetilde p_\ep\,{\rm div}_\ep(\widetilde  v_\ep)\,dz=\int_{\Omega}f'\cdot \widetilde  v_\ep'\,dz-\int_{\Omega} (\widetilde u_\ep\cdot \nabla_\ep)\widetilde u_\ep\cdot \widetilde  v_\ep\,dz,
 \end{array}
$$
which is equivalent to
$$
\begin{array}{l}
\displaystyle
(\mathcal{A}_\ep(\widetilde \varphi)-\mathcal{A}_\ep( \widetilde u_\ep),\widetilde v_\ep)-(\mathcal{A}_\ep(\widetilde \varphi),\widetilde v_\ep) +\int_{\Omega}\widetilde p_\ep\,{\rm div}_\ep(\widetilde  v_\ep)\,dz=-\int_{\Omega}f'\cdot \widetilde v_\ep'\,dz+\int_{\Omega} (\widetilde u_\ep\cdot \nabla_\ep)\widetilde u_\ep\cdot \widetilde  v_\ep\,dz.
 \end{array}
$$
Due to (\ref{Amonotone}), we can deduce 
$$
\begin{array}{l}
\displaystyle
(\mathcal{A}_\ep(\widetilde \varphi),\widetilde v_\ep) -\int_{\Omega}\widetilde p_\ep\,{\rm div}_\ep(\widetilde v_\ep)\,dz\geq \int_{\Omega}f'\cdot \widetilde  v_\ep'\,dz-\int_{\Omega} (\widetilde u_\ep\cdot \nabla_\ep)\widetilde u_\ep\cdot \widetilde  v_\ep\,dz,
 \end{array}
$$
i.e. using the expression of $\mathcal{A}_\ep$, we have
\begin{equation}\label{proof_123}
\begin{array}{l}
\displaystyle
\nu\int_{\Omega}D_\ep \widetilde \varphi:D_\ep \widetilde v_\ep\,dz+\ep^{\gamma-1}\int_{\Gamma_0}|K\widetilde \varphi'|^{s-2}K\widetilde \varphi'\cdot K\widetilde v_\ep' \,d\sigma -\int_{\Omega}\widetilde p_\ep\,{\rm div}_{\ep}(\widetilde v_\ep)\,dz\\
\noame
\displaystyle \geq \int_{\Omega}f'\cdot \widetilde  v_\ep'\,dz-\int_{\Omega} (\widetilde u_\ep\cdot \nabla_\ep)\widetilde u_\ep\cdot \widetilde  v_\ep\,dz.
 \end{array}
\end{equation}
Since $\widetilde\varphi_3\equiv 0$ and ${\rm div}_\ep(\widetilde u_\ep)=0$ in $\Omega$, it holds
$$\int_{\Omega}\widetilde p_\ep\,{\rm div}_{\ep}(\widetilde v_\ep)\,dz=\int_{\Omega}\widetilde p_\ep\,{\rm div}_{z'}(\widetilde \varphi')\,dz,$$
and from $\int_\Omega (\widetilde u_\ep\cdot \nabla_\ep \widetilde u_\ep)\widetilde u_\ep\cdot \widetilde u_\ep\,dz=0$, 
we deduce that (\ref{proof_123}) reads 
$$
\begin{array}{l}
\displaystyle
\nu\int_{\Omega}D_\ep \widetilde \varphi:D_\ep (\widetilde \varphi- \widetilde u_\ep)\,dz+\ep^{\gamma-1}\int_{\Gamma_0}|K\widetilde \varphi'|^{s-2}K\widetilde \varphi'\cdot K(\widetilde \varphi'- \widetilde u_\ep') \,d\sigma -\int_{\Omega}\widetilde p_\ep\,{\rm div}_{z'}(\widetilde \varphi')\,dz\\
\noame
\displaystyle \geq \int_{\Omega}f'\cdot (\widetilde \varphi'- \widetilde u_\ep')\,dz-\int_{\Omega} (\widetilde u_\ep\cdot \nabla_\ep)\widetilde u_\ep\cdot \widetilde\varphi\,dz.
 \end{array}
$$
Replacing $\varphi$ by $\ep^2\varphi$ and dividing by $\ep^2$ gives
$$
\begin{array}{l}
\displaystyle
\nu\int_{\Omega}\ep^2 D_\ep \widetilde \varphi:D_\ep(\widetilde \varphi- \ep^{-2}\widetilde u_\ep)\,dz+\ep^{\gamma-\gamma_s^*}\int_{\Gamma_0}|K\widetilde \varphi'|^{s-2}K\widetilde \varphi'\cdot K(\widetilde \varphi'-\ep^{-2}\widetilde u_\ep') \,d\sigma  -\int_{\Omega}\widetilde p_\ep\,{\rm div}_{z'}(\widetilde \varphi')\,dz\\
\noame
\displaystyle 
\geq \int_{\Omega}f'\cdot (\widetilde \varphi'-\ep^{-2}\widetilde u_\ep')\,dz-\int_{\Omega} (\widetilde u_\ep\cdot \nabla_\ep)\widetilde u_\ep\cdot \widetilde\varphi\,dz.
 \end{array}
$$
Next, we pass to the limit when $\ep$ tends to zero:
\begin{itemize}
\item[--] First term in the left-hand side. From convergence (\ref{velconv_power}), we get
$$\begin{array}{rl}\displaystyle 
\nu\int_{\Omega}\ep^2 D_\ep \widetilde \varphi:D_\ep(\widetilde \varphi- \ep^{-2}\widetilde u_\ep)\,dz=&\displaystyle \nu\int_{\Omega}
\partial_{z_3}\widetilde\varphi'\cdot\partial_{z_3}(\widetilde \varphi- \ep^{-2}\widetilde u_\ep')\,dz+O_\ep\\
\noame
=&\displaystyle
\nu \int_\Omega \partial_{z_3}\widetilde\varphi'\cdot\partial_{z_3}(\widetilde \varphi-  \widetilde u')\,dz+O_\ep.
\end{array}$$
\item[--] Second term in the left-hand side. Since $\gamma=\gamma_s^*$, from the continuous embedding of $H^1(0,h(z');$ $L^2(\omega))$ into $L^2(\Gamma_0)$ and convergence (\ref{velconv_power}), we get
$$\begin{array}{rl}\displaystyle \ep^{\gamma-\gamma_s^*}\int_{\Gamma_0}|K\widetilde \varphi'|^{s-2}K\widetilde \varphi'\cdot K(\widetilde \varphi'-\ep^{-2}\widetilde u_\ep') \,d\sigma=& \displaystyle \int_{\Gamma_0}|K\widetilde \varphi'|^{s-2}K\widetilde \varphi'\cdot K(\widetilde \varphi'- \ep^{-2}\widetilde u'_\ep) \,d\sigma+O_\ep\\
\noame
=& 
\displaystyle \int_{\Gamma_0}|K\widetilde \varphi'|^{s-2}K\widetilde \varphi'\cdot K(\widetilde \varphi'- \widetilde u') \,d\sigma+O_\ep.
\end{array}$$
\item[--] Third term in the left-hand side. From convergence (\ref{pressconv_power}), we get
$$\int_{\Omega}\widetilde p_\ep\,{\rm div}_{z'}(\widetilde \varphi')\,dz=\int_{\Omega}\widetilde p\,{\rm div}_{z'}(\widetilde \varphi')\,dz+O_\ep,$$
Moreover, since $\widetilde p$ is independent of $z_3$ and from condition (\ref{velconvdiv_power}), we have 
$$\int_{\Omega}\widetilde p(z')\,{\rm div}_{z'}(\widetilde u')\,dz=\int_{\omega}\widetilde p(z')\,{\rm div}_{z'}\left(\int_{0}^{h(z')}\widetilde u'\,dz_3\right)dz'=0,$$
so the third term is written as follows
$$\int_{\Omega}\widetilde p_\ep\,{\rm div}_{z'}(\widetilde \varphi')\,dz=\int_{\Omega}\widetilde p\,{\rm div}_{z'}(\widetilde \varphi'-\widetilde u')\,dz+O_\ep.$$
\item[--] First term in the right-hand side. From convergence (\ref{velconv_power}), we have
$$\int_{\Omega}f'\cdot (\widetilde \varphi'-\ep^{-2}\widetilde u_\ep')\,dz=\int_{\Omega}f'\cdot (\widetilde \varphi'- \widetilde u')\,dz+O_\ep.$$
\item[--] Second term in the right-hand side. From H${\rm \ddot{o}}$lder's inequality and estimates (\ref{velocityestimtilde_power}), we get
$$\left|\int_{\Omega}(\widetilde u_\ep \cdot \nabla_\ep)\widetilde u_\ep \widetilde \varphi\,dz\right|\leq \|\widetilde u_\ep \|_{L^2(\Omega)^3}\|D_\ep \widetilde u_\ep \|_{L^2(\Omega)^{3\times 3}}\|\widetilde \varphi\|_{L^\infty(\Omega)^3}\leq C\ep^{3},
$$
which implies
$$\int_{\Omega}(\widetilde u_\ep \cdot \nabla_\ep)\widetilde u_\ep \widetilde \varphi\,dz\to 0.$$
\end{itemize}
Finally, from previous convergences, we deduce the following limit variational inequality
$$
\begin{array}{l}
\displaystyle
\nu \int_\Omega \partial_{z_3}\widetilde\varphi'\cdot\partial_{z_3}(\widetilde \varphi-  \widetilde u')\,dz+ \int_{\Gamma_0}|K\widetilde \varphi'|^{s-2}K\widetilde \varphi'\cdot K(\widetilde \varphi'- \widetilde u') \,d\sigma - \int_{\Omega}\widetilde p\,{\rm div}_{z'}(\widetilde \varphi'-\widetilde u')\,dz
\\
\noame
\displaystyle \geq \int_{\Omega}f'\cdot (\widetilde \varphi'- \widetilde u')\,dz.
 \end{array}
$$
Since $\widetilde \varphi'$ is arbitrary, by Minty's lemma, see \cite[Chapter 3, Lemma 1.2]{Lions}, we deduce
\begin{equation}\label{limit_form_var_K}
\begin{array}{l}
\displaystyle
\nu \int_\Omega \partial_{z_3}\widetilde u'\cdot\partial_{z_3}\widetilde \varphi'\,dz+ \int_{\Gamma_0}|K\widetilde u'|^{s-2}K\widetilde u'\cdot K\widetilde \varphi' \,d\sigma - \int_{\Omega}\widetilde p\,{\rm div}_{x'}(\widetilde \varphi')\,dz
= \int_{\Omega}f'\cdot \widetilde \varphi'\,dz,
 \end{array}
\end{equation}
for every $\widetilde \varphi'\in C^1_c(\omega\times (0,h(z')))^2$ such that $\widetilde \varphi'=0$ on $\partial\Omega\setminus\Gamma_0$.  By density, this equality holds true for every $\widetilde\varphi'\in V_{z_3}$ such that $\widetilde \varphi'=0$ on $\partial\Omega\setminus\Gamma_0$. From \cite[Theorem 2.1]{Djoko},  problem (\ref{limit_form_var_K}) has a unique solution  $(\widetilde u', \widetilde p)$ and is equivalent to (\ref{limit_crit_rig})$_1$ with boundary condition (\ref{crit_bc_limit}). 
Uniqueness of solution for (\ref{limit_form_var_K}) implies that limit does not depend on the subsequence.


\end{proof}

\begin{remark} In the case $s=2$ and $K_\ep=\ep^{\gamma\over 2}\lambda^{1\over 2} I$, with $\lambda>0$, where the power slip condition (\ref{power_slip_ep}) reduces to the following Navier slip condition, with friction parameter depending on $\varepsilon$,
\begin{equation}\label{Navier_slip}-\nu\partial_{x_3} u_\ep'=-\ep^\gamma u_\ep' ,\quad  u_{\ep,3}=0,\quad  \hbox{on }   \Gamma_0.
\end{equation}
Repeating the classical proof of the existence of solution of  the Navier-Stokes problem with homogeneous Dirichlet conditions (see for instance   \cite[Theorem 2.3]{Acevedo}, \cite[Theorem 7.1]{Lions} and \cite[Theorem 10.1]{Temam}) gives the existence of at least a weak solution $(u_\ep, p_\ep)\in \mathbb{V}(\Omega^\ep)\times L^2_0(\Omega^\ep)$  of problem (\ref{system_1introad}) and (\ref{Navier_slip}).
\\

Proceeding similarly to the proof of Theorem \ref{thm1_critical_rig}, but without monotonicity arguments, we can prove that  the limit system satisfied by $(\widetilde u', \widetilde p)$ is (\ref{limit_crit_rig}), and that there exists a critical value of $\gamma$ is $-1$, which agrees with $\gamma_2^*$, such that the effective limit conditions on $\Gamma_0$ are the following ones
\begin{itemize}
\item If $\gamma=-1$, then 
$$-\nu\partial_{z_3}\widetilde u'=-\lambda \widetilde u'\quad\hbox{on }\Gamma_0.$$
\item If $\gamma>-1$, then 
$$-\nu\partial_{z_3}\widetilde u'=0\quad\hbox{on }\Gamma_0.$$
\item If $\gamma<-1$, then 
$$\widetilde u'=0\quad\hbox{on }\Gamma_0.$$
\end{itemize}
By classical arguments of the existence and uniqueness of solution for the Stokes problem with homogeneous Dirichlet conditions, there exists a unique weak solution $(\widetilde u', \widetilde p)$ of problem (\ref{limit_crit_rig}) with corres\-ponding boundary conditions on $\Gamma_0$ depending on the value of $\gamma$ given above.
\end{remark}

\section*{Acknowledgments}
 We would like to thank all those people (publishers, editors, referees  and researchers) for all the support received for the deve\-lopment of our lines of research. In particular, Mar\'ia would like to dedicate this article to her father, Julio, for his unconditional support.

 \section*{Conflict of interest}
 The authors confirm that there is no conflict of interest to report.

 \section*{Data availability statement}
Data sharing not applicable to this article as no datasets were generated or analysed during the current study.


\begin{thebibliography}{00}


\bibitem{Acevedo} P. Acevedo, C. Amrouche, C. Conca, and A. Ghosh, {\it Stokes and Navier--Stokes equations with Navier boundary conditions}, J. Differ. Equ. {\it 285} (2021), 258--320

\bibitem{Aldbaissy} R. Aldbaissy, N. Chalhoub, J.K. Djoko, and T. Sayah, {\it Full discretisation of the time dependent Navier-Stokes equations with anisotropic slip boundary condition}, Int. J. Numer. Anal. Mod. {\bf 20} (2023),  497--517.
%
\bibitem{Aldbaissy2} R. Aldbaissy, N. Chalhoub, J.K. Djoko, and T. Sayah, {\it Full discretization of the time dependent Navier-Stokes equations with anisotropic slip boundary condition coupled with the convection-diffusion-reaction equation}, SeMA (2024).

\bibitem{Amrouche} C. Amrouche and A. Rejaiba, {\it $L^p$-theory for Stokes and Navier--Stokes equations with Navier boundary condition}, J. Differ. Equ. {\bf 256} (2014), 1515--1547.


%
\bibitem{AnguianoSGNonlinear} M. Anguiano and F.J. Suárez-Grau,  {\it Nonlinear Reynolds equations for non-Newtonian thin-film fluid flows over a rough boundary}, IMA J. Appl. Math. {\bf 84}  (2019), 63--95.


\bibitem{AnguianoSGmag} M. Anguiano and F.J. Suárez-Grau, {\it Mathematical derivation of a Reynolds equation for magneto-micropolar fluid flows through a thin domain} {\bf 75} (2024), 75: 28.
%
\bibitem{Bayada} G. Bayada and M. Chambat,  {\it The transition between the Stokes equations and the Reynolds equation: A mathematical proof}, Appl. Math. Optim.  {\bf 14} (1986), 73--93.
%
 \bibitem{Bonnivard} M. Bonnivard, A.-L. Dalibard, and D. G\'erard-Varet, {\it Computation of the effective slip of rough hydrophobic surfaces via homogenization}, Math. Mod.  Meth. Appl. S. {\bf 24} (2014), 2259--2285.

 \bibitem{Bonnivard_SG1} M. Bonnivard and F.J. Su\'arez-Grau, {\it On the influence of wavy riblets on the slip behaviour of viscous fluids}, Z. Angew. Math. Phys. {\bf 67} (2016), 67: 27. 

 \bibitem{Bonnivard_SG2} M. Bonnivard and F.J. Su\'arez-Grau, {\it Homogenization of a Large Eddy Simulation Model for Turbulent Fluid Motion Near a Rough Wall}, J. Math. Fluid Mech. {\bf 20} (2018), 1771--1813. 
 


\bibitem{Boyer} F. Boyer and P. Fabrie, {\it Mathematical Tools for the Study of the Incompressible Navier-Stokes Equations and Related Models},  Springer Science \& Business Media, 2013.


\bibitem{BucurNN} D. Bucur, E. Feireisl, and S. Ne${\rm \check{c}}$asov\'a, {\it Influence of wall roughness on the slip behaviour of viscous fluids}, Proc. R. Soc. Edinb. A: Math. {\bf 138} (2008), 957--973.


\bibitem{Bucur} D. Bucur, E. Feireisl,  S. Ne${\rm \check{c}}$asov\'a, and J. Wolf, {\it On the asymptotic limit of the Navier--Stokes system on domains with rough boundaries}, J. Differ. Equ. {\bf 244} (2008), 2890--2908.




\bibitem{Clupeau} T. Clopeau, A. Mikeli\'c, and R. Robert, {\it On the vanishing viscosity limit for the 2D incompressible Navier--Stokes equations with the friction type boundary conditions}, Nonlinearity  {\bf 11} (1998), 1625--1636.

\bibitem{Dalibard} A.-L. Dalibard and D. G\'erard-Varet, {\it Effective boundary condition at a rough surface starting from a slip condition}, J. Differ. Equ.  {\bf 251} (2011), 3297--3658.

\bibitem{Djoko} J.K. Djoko, J. Koko, M. Mbehou, and T. Sayah, {\it Stokes and Navier-Stokes equations under power law slip boundary condition: Numerical analysis}, Comput. Math. Appl. {\bf 128} (2022), 198--213.


\bibitem{Djoko2} J. K. Djoko, V. S. Konlack, and T. Sayah, {\it Power law slip boundary condition for Navier-Stokes equations: Discontinuous Galerkin schemes}, Comput. Geosci. {\bf 28} (2024), 107--127.

\bibitem{Galdi} G.P. Galdi,  {\it An Introduction to the Mathematical Theory of the Navier-Stokes Equations},  Springer,
New York, 1994.
%
\bibitem{Lions} J.L. Lions, {\it Quelques m\'ethodes de r\'esolution des probl\`emes aux limites non lin\'eaires}, Dunod,
Gauthier-Villars, Paris, 1969.


\bibitem{Mikelic} A. Mikeli\'c, {\it An Introduction to the Homogenization Modeling of Non-Newtonian and Electrokinetic Flows in Porous Media}, A. Farina, A. Mikeli\'c, F. Rosso, eds, Non-Newtonian Fluid Mechanics and Complex Flows. Lecture Notes in Mathematics, vol 2212, Springer,  2018, pp 171--227.
%
\bibitem{MikTap} A. Mikeli\'c  and R. Tapiero,  {\it Mathematical derivation of the power law describing polymer flow through a thin slab}, RAIRO Mod\'el. Math. Anal. Num\'er. {\bf 29} (1995), 3--21.

 \bibitem{Navier} C.L.M.H. Navier, {\it Sur les lois d'\'equilibre et du mouvement des corps \'elastiques}, M\'em. Acad. Sci.  {\bf 7} (1827), 375--394.
%

\bibitem{Solonnikov} V.A. Solonnikov and V.E. ${\rm \check{S}}$${\rm \check{c}}$adilov, {\it A certain boundary value problem for the stationary system of Navier--Stokes equations}, Tr. Mat. Inst. Steklova {\bf 125} (1973), 196--210.
%
\bibitem{SG1} F.J. Su\'arez-Grau,  {\it Effective boundary condition for a quasi-newtonian fluid at a slightly rough boundary starting from a Navier condition}, ZAMM--Z. Angew.  Math. Me. {\bf 95} (2015), 527--548.

 \bibitem{SG2}F.J. Su\'arez-Grau, {\it Asymptotic behavior of a non-Newtonian flow in a thin domain with Navier law on a rough boundary}, Nonlinear Anal-Theor. {\bf 117} (2015), 99--123.


\bibitem{Temam} R. Temam, {\it Navier-Stokes Equations}, North Holland, 1984.


\end{thebibliography}
\end{document}